\documentclass[number,citesort,dvips]{arxbj}
\usepackage{upgreek}
\usepackage{graphicx}

\aid{0}
\volume{17}
\issue{1}
\pubyear{2011}
\firstpage{290}
\lastpage{319}
\doi{10.3150/10-BEJ269}

\makeatletter

\newremark{remark}{Remark}[section]

\newtheorem{corollary}{Corollary}[section]
\newtheorem{theorem}{Theorem}[section]
\newtheorem{lemma}{Lemma}[section]

\newproclaim{definition}{Definition}[section]

\newcommand{\point}{(\nu,\omega) \in(0, 2\uppi]^2}
\renewcommand{\Re}{\operatorname{Re}}
\renewcommand{\Im}{\operatorname{Im}}

\newcommand{\cov}{\operatorname{cov}}
\newcommand{\cosec}{\operatorname{cosec}}
\newcommand{\norm}[1]{\Vert #1  \Vert}
\newcommand{\ciag}[1]{\{#1 \dvtx t\in\mathbb{Z}\}}

\newcommand{\im}{\mathrm{i}}

\makeatother

\begin{document}
\begin{frontmatter}

\title{Asymptotic distributions and subsampling in spectral analysis
for almost periodically correlated time series}
\runtitle{Spectral analysis for APC time series}

\begin{aug}
\author{\fnms{\L{}ukasz} \snm{Lenart}\ead[label=e1]{llenart@wsb-nlu.edu.pl}\corref{}} 
\runauthor{\L{}. Lenart}
\address{Department of Econometrics,
The Graduate School of Business WSB-NLU,
ul. Zielona 27, 33-300 Nowy Sacz, Poland. \printead{e1}}
\end{aug}

\received{\smonth{11} \syear{2008}}
\revised{\smonth{11} \syear{2009}}

%
\begin{abstract}
The aim of this article is to establish asymptotic distributions
and consistency of subsampling for spectral density and for
magnitude of coherence for non-stationary, almost periodically
correlated time series. We show the asymptotic normality of the
spectral density estimator and the limiting distribution of a
magnitude of coherence statistic for all points from the
bifrequency square. The theoretical results hold
under $\alpha$-mixing and moment conditions.
\end{abstract}

%

\end{frontmatter}

\section{Introduction}\label{sec1}

The analysis of the second order structure of time series is usually
based on the characterization of the autocovariance function. One
possibility is that the time series is not stationary, with a periodic
or an almost periodic autocovariance function. The models with this
structure were applied in many fields, including telecommunications
\cite{gardner,Nap01}, meteorology
\cite{BloHurdLund94}, econometrics \cite{parzen,osborn} and many
others (see the review in \cite{giannakis05}
or \cite{gardner2}). This class of non-stationary time series
was introduced by \cite{gladyshev1}. Formally, we say that a
second-order real-valued time series
$\{X_{t}\dvtx t\in\mathbb{Z}\}$ is \textit{periodically correlated}
(PC for
short) if both the mean $\mu(t)=E(X_{t})$ and the autocovariance
function $B(t,\tau)=\operatorname{cov}(X_{t},X_{t+\tau})$ are periodic at
$t$ for every $\tau\in\mathbb{Z}$, with the same period $T$. A~second-order real-valued time series $\ciag{X_{t}}$ is called
\textit{almost periodically correlated} (APC for short) if both the mean
$\mu(t)$ and the autocovariance function $B(t,\tau)$ are almost
periodic function at $t$ for every $\tau\in\mathbb{Z}$. The
definition of almost periodic function that we use in this paper
can be found in \cite{Cor}, page 45. It is easy to see that the
class of APC time series contains the class of PC time series and
the class of stationary time series.

One of the main problems is how to detect this kind of non-stationary
structure in the time series. This problem was
considered in the time domain by Vecchia and Ballerini \cite
{vecchia91}, Dandawate and Giannakis \cite{giannakis94}, %
Dehay and Le\'skow \cite{Dehles96b}, Synowiecki \cite{Syn07} and,
recently, by Lenart \textit{et al.} \cite{Len08b}.
Lenart \textit{et al.} \cite{Len08b} present tests for PC structure
based on the estimator of
Fourier coefficient $a(\lambda,\tau)$ and the subsampling methodology.
In the frequency domain, the problem of detecting periodicity was
considered by Hurd and Gerr \cite{hurdgra}, Lund \textit{et al.} \cite
{LundHurdBlo95}, Broszkiewicz-Suwaj \textit{et al.} \cite{Brosz},
Lii and Rosenblatt \cite{Rosenblatt} and Hurd and Miamee \cite
{hurdmiamee}. Hurd and Gerr \cite{hurdgra} and
Hurd and Miamee \cite{hurdmiamee} present a graphical method for detecting
stationarity and periodicity, where the testing statistics were
approximated by a beta distribution. However, this method has not yet
been justified for non-Gaussian white noise or dependent random
variables. Broszkiewicz-Suwaj \textit{et al.} \cite{Brosz} have used
the \textit{moving block bootstrap}
to construct methods for detecting periodicity in time series. Still,
the consistency of the bootstrap
applied for their testing statistics remains an open problem. Lii and
Rosenblatt \cite{Rosenblatt} present an innovative
algorithm for estimating the support of the spectral measure for
Gaussian APC processes. The authors considered the case where
the support of the spectral measure is concentrated on a finite number
of parallel lines.

In the next section, we recall the theoretical background and formulate
the assumptions which are used in the
sections which follow. We consider general APC time series for which
the support of
the spectral measure can be concentrated on an infinite
number of parallel lines. Section~\ref{sec3} presents the consistent estimator
for the magnitude of coherence and the
extension of the spectral
density function to the bifrequency square $(0,2\uppi]^2$. In this
section, we show the asymptotic normality of
the normalized spectral density estimator, the asymptotic distribution
of the normalized magnitude of
the spectral density estimator and the coherence statistic for all
points from
the bifrequency square $(0,2\uppi]^2$. We give the exact forms for
these distributions. It is shown that the asymptotic
distribution of the normalized magnitude of the coherence statistic
strongly depends on the support of the spectral measure. This
distribution is not Gaussian for the points which do
not belong to the support of the spectral measure. Recall that for
the time series, the asymptotic normality of the spectral density
estimator was shown for the stationary case in \cite{zurbenko} under
an $\alpha$-mixing condition and in \cite{priestley} for linear
filters. For the PC case, the asymptotic normality was established in
\cite{hurdmiamee} for linear filters. Asymptotic distributions from
Section \ref{sec3} provide a possibility to establish consistency of the
subsampling for the magnitude of spectral density and coherence,
which is shown in Section \ref{sec4}. It should be emphasized that magnitude of
spectral density
and coherence are broadly used as fundamental characteristics for
telecommunication signals (see \cite{gardner86}).
Section \ref{sec5} presents a simulation study, where we show applications of
our results.
Following the ideas of Hurd and Gerr \cite{hurdgra} and Hurd and
Miamee \cite{hurdmiamee}, we present graphical methods for detecting
periodicity
in autocovariance and give theoretical justifications for these methods.
Finally, the asymptotic distribution and consistency of subsampling
provide a possibility to construct asymptotically consistent confidence
intervals for magnitude of spectral density and coherence. All proofs
can be found in the \hyperref[app]{Appendix}.

\section{Definitions and assumptions}\label{sec2}

We will be focusing on second-order inference for time series,
therefore, for simplicity, we make the following
assumption:
\begin{longlist}[(A1)]
\item[(A1)] \textit{Assume that the real-valued time series $\ciag
{X_{t}}$ is zero-mean.}
\end{longlist}

Moreover, we will consider non-stationary, almost periodically
correlated time series; therefore, we formulate the next
assumption:
\begin{longlist}[(A2)]
\item[(A2)] \textit{Assume that the time series $\ciag{X_{t}}$ is APC.}
\end{longlist}

In the APC case with $\mu(t) \equiv0$
for any $\tau\in\mathbb{Z}$, the autocovariance function $B(t,\tau)$
has the Fourier representation
\begin{equation} \label{fourier}
B(t,\tau) \sim\sum_{\lambda\in\Lambda_{\tau}}a(\lambda,\tau)\mathrm
{e}^{\mathrm{i} \lambda t},
\end{equation}
where $a(\lambda,\tau)$ are the Fourier coefficients of the form
$a(\lambda,\tau)=\lim_{n \to
\infty}n^{-1}\sum_{j=1}^{n}B(j,\tau)\mathrm{e}^{-\mathrm{i} \lambda j}$
and the set $\Lambda_{\tau}= \{\lambda\in[0,2\uppi)\dvtx a(\lambda
,\tau)\not= 0 \}$ is countable (see \cite{Cor}).

Define the set $\Lambda=\bigcup_{\tau\in\mathbb{Z}} \Lambda_{\tau}$.
We make the following assumption concerning the
set $\Lambda$ and the Fourier coefficients $a(\lambda,\tau)$.
We need this assumption in order to establish asymptotic properties of
the spectral density
estimator:
%
\begin{longlist}[(A3)]
\item[(A3)] \textit{Assume that for any $x \in[0,2 \uppi)$, there
exists a real-valued sequence
$\{z_{\tau}(x)\}_{\tau\in\mathbb{Z}}$ such that we have
\begin{equation}
\sum_{\lambda\in\Lambda_{\tau} \setminus\{x \} }
\biggl|a(\lambda,\tau) \cosec\biggl({\frac{\lambda- x}{2}} \biggr)
\biggr|
< z_{\tau}(x) < \infty
\label{fourierx}
\end{equation}
and $z_{\tau}(x) \to0 $ for $|\tau| \to\infty$.}
\end{longlist}

\begin{remark}
By (A3), we have that $\sum_{\lambda\in\Lambda_{\tau} } |a(\lambda
,\tau) | < \infty,$ which implies that
\begin{equation} \label{fourier2}
B(t,\tau) = \sum_{\lambda\in\Lambda_{\tau}}a(\lambda,\tau)\mathrm
{e}^{\mathrm{i} \lambda t}.
\end{equation}
If $\operatorname{card}(\Lambda)<\infty$, and we assume that $|B(t,\tau
)|<p_{\tau}$ uniformly at $t$ and $p_{\tau} \to0$ for
$|\tau|\to\infty$, then (A3) holds. Note that (\ref{fourierx}) implies
that for any $\lambda\in[0, 2 \uppi)$, we have that $|a(\lambda,\tau
)|\to0$ as $|\tau| \to\infty$.
\end{remark}

To introduce the spectral theory for APC time series, we need the following
assumption:

\begin{longlist}[(A4)]
\item[(A4)] \textit{Assume that the time series $\ciag{X_{t}}$ is
harmonizable (for the definition, see} \cite{loeve} \textit{or}
\cite{yaglom}\textit{)}.
\end{longlist}

By assumptions (A2) and (A4), the spectral measure $R$, defined on the
bifrequency square $(0,2 \uppi]^2$,
has support contained in the set (see \cite{Dehhurd93})
\[
S= \bigcup_{|\lambda| \in\Lambda} \{(\nu,\omega)\in(0,2 \uppi
]^2\dvtx\omega=\nu-\lambda\}.
\]

Moreover, the coefficients $a(\lambda,\tau)$ are the Fourier
transforms of complex measures $r_{\lambda}(\cdot)$, that is,
$a(\lambda,\tau)=\int_{0}^{2\uppi}\mathrm{e}^{\mathrm{i} \xi\tau
}r_{\lambda}(\mathrm{d}\xi)$. The measure $r_{\lambda}$ can be identified
with the
restriction of the spectral measure $R$ to the line $\omega=\nu
-\lambda$, where $|\lambda| \in\Lambda$. We need the following
assumption:
\begin{longlist}[(A5)]
\item[(A5)] \textit{Assume that the measure $r_{0}$ is absolutely
continuous with respect to the Lebesgue measure.}
\end{longlist}

By the above assumptions, for any $|\lambda| \in\Lambda$, there exists
a spectral density function $g_{\lambda}(\cdot)$ such that
$ g_{\lambda}(\nu)=\frac{1}{2 \uppi}\sum_{\tau=-\infty}^{\infty}a(\lambda
,\tau)\mathrm{e}^{-\mathrm{i} \nu\tau}$
(see \cite{Hurcor91} and \cite{Dehhurd93}).

Now, let us define the extension of the spectral density function
$P(\cdot,\cdot)$ to the bifrequency square $(0, 2 \uppi]^2$ via
\begin{equation}\label{as-per}
P(\nu,\omega)= \cases{
g_{\lambda}(\nu) &\quad  for $\omega=\nu- \lambda $ and $|\lambda| \in\Lambda$, \cr
0 &\quad  for $(\nu,\omega)\notin S$.
}
\end{equation}
Hence,
\begin{equation} \label{spec_rec}
P(\nu,\omega) = \frac{1}{2
\uppi}\sum_{\tau=-\infty}^{\infty}a(\nu-\omega,\tau)\mathrm{e}^{-\mathrm{i} \nu\tau}
\end{equation}
for any $\point$, which follows from the fact that $a(\nu-\omega,\tau
)=0$ for any $|\nu-\omega| \in[0,2\uppi) \setminus\Lambda_{\tau}$
and $\tau\in\mathbb{Z}$. If we assume that $g_{0}(\xi)\not= 0$ for
any $\xi \in(0, 2 \uppi]$, then for any point $(\nu,\omega)$ from the
bifrequency square $ (0,2\uppi]^2$, we define a magnitude of coherence via
\begin{equation}\label{cohe}
|\gamma(\nu,\omega)|= \cases{\displaystyle
\frac{|g_{\lambda}(\nu)|}{\sqrt{g_{0}(\nu)g_{0}(\omega)}} & \quad for $ \omega=\nu- \lambda $ and $ |\lambda| \in\Lambda$,\cr
0 & \quad for $ (\nu,\omega) \notin S$.
}
\end{equation}

Note that the magnitude of coherence is a real number from the interval
$[0,1]$, while the extension of the spectral density function is a
complex number.

\section{Asymptotic distributions}\label{sec3}

The estimator of $P(\nu,\omega)$ (based on a sample $\{
X_{1},X_{2},\ldots,X_{n}\}$) of the form
\begin{equation} \label{g-r}
\hat P_{n}(\nu,\omega)=\frac{1}{2\uppi n}\sum_{s=1}^{n}\sum
_{t=1}^{n}X_{s}X_{t}\mathrm{e}^{-\mathrm{i} \nu s}\mathrm{e}^{\mathrm
{i} \omega t}
\end{equation}
is generally not consistent in the mean square sense for APC time
series. The exact form of the asymptotic
variance of this estimator was calculated in \cite{Len08a}.

To obtain a consistent
estimator, we consider the following class of smoothed
estimators of $P(\nu,\omega)$ (see \cite{zurbenko} for the stationary case):
\begin{equation} \label{dens}
\hat G_{n}(\nu,\omega)=\frac{1}{2\uppi n}\sum_{s=1}^{n}\sum_{t=1}^{n}
H_{L_{n}}(s-t)X_{s}X_{t}\mathrm{e}^{-\mathrm{i} \nu s}\mathrm
{e}^{\mathrm{i} \omega t},
\end{equation}
where $H_{L_{n}}(\cdot)$ is a lag window function such that
$H_{L_{n}}(\tau)=0$ for $|\tau|>L_{n}$, $\tau\in\mathbb{Z}$ and
$L_{n}$ is
a sequence of positive integers tending to infinity with $n$. Moreover,
we assume that $L_{n}/n\rightarrow0$.
If we assume that $g_{0}(\xi)\not= 0$ for any $\xi\in(0,
2\uppi]$, then we define a magnitude of coherence statistic based
on the estimator (\ref{dens}) via
\begin{equation} \label{cohth}
|\hat\gamma_{n}(\nu,\omega)|=\frac{|\hat G_{n}(\nu,\omega)|}{\sqrt
{|\operatorname{Re}(\hat G_{n}(\nu,\nu))\operatorname{Re}(\hat
G_{n}(\omega,\omega))|}}.
\end{equation}

The following considerations will involve the general assumption that
we have a sample
$\{X_{c_{n}+1},X_{c_{n}+2},\ldots,X_{c_{n}+d_n}\}$ from time series
$\ciag{X_{t}}$, where $\{c_{n}\}_{n \in
\mathbb{Z}}$ and $\{d_{n}\}_{n \in\mathbb{Z}}$ are arbitrary sequences
of non-negative integers such that $d_n \to\infty$. This
will help us to prove the consistency of the subsampling procedure to
be discussed in Section \ref{sec4}.

Given a sample $\{X_{c_{n}+1},X_{c_{n}+2},\ldots,X_{c_{n}+d_n}\}$, we
define the estimators of $P(\nu,\omega)$
and $|\gamma(\nu,\omega)|$ via
\begin{eqnarray}
\hat P_{n}^{c,d}(\nu, \omega)&=&\frac{1}{2\uppi d_n}\sum_{s=c_n
+1}^{c_n +d_n}\sum_{t=c_n +1}^{c_n +d_n} X_{s}X_{t}\mathrm{e}^{-\mathrm
{i} s \nu}\mathrm{e}^{\mathrm{i}
t \omega}, \label{den-p}
\\
\hat G_{n}^{c,d}(\nu, \omega)&=&\frac{1}{2 \uppi d_n}\sum_{s=c_n
+1}^{c_n +d_n}\sum_{t=c_n +1}^{c_n +d_n}
H_{L_{d_n}}(s-t)X_{s}X_{t}\mathrm{e}^{-\mathrm{i} s \nu}\mathrm
{e}^{\mathrm{i} t \omega}
\label{den-h}
\end{eqnarray}
and
\begin{equation}\label{cohthphi}
|\hat\gamma_{n}^{c,d}(\nu,\omega)|=\frac{|\hat G_{n}^{c,d}(\nu,\omega
)|}{\sqrt{|\operatorname{Re}(\hat G_{n}^{c,d}(\nu,\nu))\operatorname
{Re}(\hat G_{n}^{c,d} (
\omega,\omega))|}}.
\end{equation}
Besides the fact that, in this work, we consider zero-mean time
series, in the following remark, we define the estimator of the
spectral density
function in the case where the mean function $\mu(t)$ is an almost
periodic function such that $\mu(t)\not\equiv0$.

\begin{remark}
Assume that the mean function
$\mu(t)=E(X_t)$ for APC time series is
an almost periodic function such that $\mu(t)=\sum_{\gamma\in\Gamma}
b(\gamma) \mathrm{e}^{\mathrm{i} \gamma t}$, where $b(\gamma)$ are the
Fourier coefficients and the
set $\Gamma$ is known and finite. The natural estimator of the mean function
$\mu(t)$ (for $c_n+1 \leq t \leq c_n+d_n$) based on a sample $\{
X_{c_n+1},X_{c_n+2},\ldots,X_{c_n+d_n} \}$ then takes the form
$\hat\mu_{n}^{c,d}(t) = \sum_{\gamma\in\Gamma} \hat
b_{n}^{c,d}(\gamma) \mathrm{e}^{\mathrm{i} \gamma t}$,
where $\hat b_{n}^{c,d}(\gamma)=\frac{1}{d_n}\sum
_{j=c_n+1}^{c_n+d_n}X_{j}\mathrm{e}^{-\mathrm{i} \gamma j}$
is an estimator of the Fourier coefficient $b(\gamma)$ for any $\gamma
\in\Gamma$. The more general natural
estimator of the spectral density
$P(\nu,\omega)$ based on a sample $\{ X_{c_n+1},X_{c_n+2},\ldots
,X_{c_n+d_n} \}$ then takes the form
\begin{equation}\label{dens2}
\hat T_{n}^{c,d}(\nu,\omega)=\frac{1}{2\uppi d_n}\sum
_{s=c_n+1}^{c_n+d_n}\sum_{t=c_n+1}^{c_n+d_n}
H_{L_{d_n}}(s-t)\bigl(X_{s}-\hat\mu_{n}^{c,d}(s)\bigr)\bigl(X_{t}-\hat\mu
_{n}^{c,d}(t)\bigr)\mathrm{e}^{-\mathrm{i} \nu s}\mathrm{e}^{\mathrm{i}
\omega t}.
\end{equation}
Similar modification in the PC case can be found in
\cite{hurdmiamee}, Section 10.4. \label{non-zero}
\end{remark}

In Theorem \ref{covariance}, we show the
asymptotic covariance between normalizing estimators $\hat
G_{n}^{c,d}(\nu_1 , \omega_1 )$ and $\hat G_{n}^{c,d}(\nu_1 , \omega_1
)$ for
$(\nu_1,\omega_1)$, $(\nu_1,\omega_1) \in[0,2\uppi)^2$. We need the
following assumption concerning the lag window function $H_{L_{n}}(\cdot)$.

\begin{longlist}[(B)]
\item[(B)] \textit{Assume that for the lag window function $H_{L_n}(\cdot
)$, where $L_{n}\to\infty$, $L_{n}/n \to0$, there exists a
real-valued function $w(\cdot)$ such that:
\begin{longlist}[(iii)]
\item[(i)] $w(x)=0$ for $|x|>1$ and $w(\cdot)$ is an even
function and non-increasing on the interval $[0,1]$;
\item[(ii)] there exists a real number $\theta\in(0,1]$
such that for any $|x|\leq\theta$, we have $w(x)=1$;
\item[(iii)] $H_{L_{n}}(\tau)=w(\tau/L_n)$ for any $\tau
\in\mathbb{Z}$ and $n \geq1$;
\item[(iv)] $w(\cdot)$ is
Lipschitz continuous on the interval $[-1,1]$ with Euclidean
metric, which means that there exists a real number $W$ such that
for any $x,y \in[-1,1]$, we have
\[
|w(x)-w(y)|\leq W |x-y|.
\]
\end{longlist}
}
\end{longlist}

For the convenience of the reader, before Theorem \ref{covariance}, we introduce the
concept of $\alpha$-mixing.

\begin{definition}[(\cite{doukhan})]
Let $\ciag{X_t}$ be a real-valued time series. The $\alpha$-mixing
sequence $\alpha(\cdot)$ which corresponds to the time series $\ciag
{X_t}$ is
defined as
\[
\alpha(s)=\sup_{t\in\mathbb{Z}} \mathop{\sup_{A\in\mathcal
{F}_X(-\infty,
t)}}_{ B\in\mathcal{F}_X(t+s, \infty)} |P(A\cap
B)-P(A)P(B)|,
\]
where $\mathcal{F}_X(t_1, t_2)$ stands for the
$\sigma$-algebra generated by $\{X(t) \dvtx t_1\leq t\leq t_2\}$. The time
series $\ciag{X_t}$ is called $\alpha$-mixing if
$\alpha(s) \to0$ as $s \to\infty$.
\end{definition}
For any random variable $X$ and positive constant $p$, we define the norm
$\|X\|_{p}=(E|X|^p)^{1/p}$.

\begin{theorem}\label{covariance}
Assume that \textup{(A1)--(A5)} and \textup{(B)} hold. If, additionally, there
exist $\delta>0$, $\Delta< \infty$ and $K < \infty$ such that:
\begin{itemize}[(ii)]
\item[(i)]
$\sup_{t \in\mathbb{Z}} \norm{X_{t}}_{6+3\delta} < \Delta;$\vspace*{2pt}
\item[(ii)]
$\sum_{k=0}^{\infty}(k+1)^2\alpha(k)^{\delta/{(2+\delta)}}\leq K$,
\end{itemize}
then for any $(\nu_{1},\omega_{1}),(\nu_{2},\omega_{2})\in(0,2\uppi]^2$,
we have the convergence
\begin{eqnarray*}
&&\lim_{n \to\infty} \frac{d_n}{L_{d_n}} \operatorname{cov} (\hat
G_{n}^{c,d}(\nu_{1} ,\omega_{1}),\hat G_{n}^{c,d}(\nu_{2}
,\omega_{2}) ) \\
&&\quad = \rho\bigl( P(\nu_{1},\nu_{2}) \overline{P(\omega_{1},\omega
_{2})} + P(\nu_{1},2\uppi - \omega_{2}) \overline{P(\nu_{2},2\uppi
-\omega_{1})} \bigr),
\end{eqnarray*}
where $\rho=\int_{-1}^{1}w^2(x) \,\mathrm{d}x$.
\end{theorem}

In Theorem \ref{normality}, we show that the estimator
(\ref{den-h}) has an asymptotically normal distribution. This result
is crucial for establishing the asymptotic
distribution for the magnitude of coherence statistic and proving the
consistency of subsampling applied
for magnitude of coherence and spectral density.

\begin{theorem}\label{normality}
Assume that \textup{(A1)--(A5)} and \textup{(B)} hold. Additionally,
assume that:
\begin{longlist}[(iii)]
\item[(i)] there exists $\delta>0$ such that $\sup_{t
\in\mathbb{Z}} \norm{X_{t}}_{6+3\delta}\leq\Delta<\infty;$
\item[(ii)] $L_{n}=\mathrm{O}(n^{\kappa})$ for some $\kappa\in
(0,\delta/(4+4\delta))$;
\item[(iii)]
$\sum_{k=0}^{\infty}(k+1)^{r}\alpha(k)^{\delta/{(r+2+\delta)
}}<K<\infty$, where $r$ is the even integer such that
\[
r>\max \biggl\{1+3 \delta/2, \frac{1-\kappa}{ 2 \kappa },\frac{2 \kappa(1+\delta)}{\delta-2 \kappa(1+\delta)} \biggr\}.
\]
\end{longlist}
Then
\begin{equation} \label{eq32}
\sqrt{\frac{d_n}{L_{d_n}}} \left( \left[
\matrix{
\operatorname{Re}(\hat G_{n}^{c,d}(\nu,\omega))\cr
\operatorname{Im}(\hat G_{n}^{c,d}(\nu,\omega))
}\right] - \left[
\matrix{
\operatorname{Re}(P(\nu,\omega))\cr
\operatorname{Im}(P(\nu,\omega))
}
\right] \right) \stackrel{d}{\longrightarrow} \mathcal{N}_{2}(0,\Sigma(\nu,\omega)),
\end{equation}
where $\Sigma(\nu,\omega)=[\sigma_{ij}]_{i,j=1,2}$,
\begin{eqnarray}
\sigma_{11} & =& \tfrac{1}{2} \bigl( g_{0}(\nu)g_{0}(\omega) + |P(\nu,2 \uppi-
\omega)|^2 +
\operatorname{Re} \bigl( P(\nu,2 \uppi- \nu) P(2\uppi- \omega, \omega) \bigr)\nonumber\\
&&\hphantom{\frac{1}{2} \bigl(}{}
+ [\operatorname{Re}(P(\nu,\omega))]^2- [\operatorname{Im}(P(\nu,\omega))]^2 \bigr),\nonumber\\
\sigma_{12} & =& \sigma_{21} = - [\operatorname{Re}(P(\nu,\omega))]^2 [
\operatorname{Im}(P(\nu,\omega))]^2
- \tfrac{1}{2}\operatorname{Im} \bigl(P(\nu,2\uppi- \nu)P(2\uppi- \omega,\omega) \bigr),\qquad \\
\sigma_{22} & =& \tfrac{1}{2} \bigl( g_{0}(\nu)g_{0}(\omega) + |P(\nu,2 \uppi-
\omega)|^2 - \operatorname{Re} \bigl( P(\nu,2 \uppi- \nu) P(2\uppi- \omega, \omega) \bigr)\nonumber\\
&&\hphantom{\frac{1}{2} \bigl(}{}
- [\operatorname{Re}(P(\nu,\omega))]^2+ [\operatorname{Im}(P(\nu,\omega))]^2 \bigr).\nonumber
\end{eqnarray}
%
\end{theorem}

Recall that in the stationary case, for the majority of spectral
windows, the optimal $L_n$ (in
the mean square sense) for the estimator $\hat G_{n}(\nu,\omega)$ is of
order $n^{1/5}$
(see \cite{priestley}, pages 462--463, or \cite{Pol}, page 86). Note that
condition (ii) of Theorem \ref{normality} includes the case
$L_{n}=\mathrm{O}(n^{1/5})$ for $\delta> 4$.

The following corollary is a natural generalization of the previous
theorem to the multidimensional case. We need this result to establish the
asymptotic distribution for the magnitude of coherence statistic.

\begin{corollary}
Let $\ciag{X_{t}}$ be a time series such that all the assumptions of
Theorem \ref{normality} hold.
We then have the convergence
\begin{equation}\label{gamma1}
\sqrt{\frac{d_n}{L_{d_n}}} \left( \left[
\matrix{
\operatorname{Re}(\hat G_{n}^{c,d}(\nu,\omega))\cr
\operatorname{Re}(\hat G_{n}^{c,d}(\nu,\nu))\cr
\operatorname{Re}(\hat G_{n}^{c,d}(\omega,\omega))\cr
\operatorname{Im}(\hat G_{n}^{c,d}(\nu,\omega))
}
\right] - \left[
\matrix{
\operatorname{Re}(P(\nu,\omega))\cr
g_{0}(\nu) \cr
g_{0}(\omega)\cr
\operatorname{Im}(P(\nu,\omega))
}
\right] \right) \stackrel{d}{\longrightarrow} \mathcal{N}_{4}(0,{\Psi}(\nu,\omega
)),
\end{equation}
where the covariance matrix $\Psi(\nu,\omega)$ can be obtained by Theorem
\ref{covariance}. \label{mnormality}
\end{corollary}

\begin{remark}
Note that in the special case where $\nu=\omega$, we have $
\operatorname{rank}(\Sigma(\nu,\nu))=1$ and $\operatorname{rank}(\Psi(\nu,\nu))=3$, which
follow from the fact that
$\operatorname{Im}(\hat G_{n}^{c,d}(\nu,\nu))=0$.
\end{remark}

Let us establish the asymptotic distribution of the normalized
estimator $|\hat G_{n}^{c,d}(\nu,\omega)|$
(under all the assumptions of
Theorem \ref{normality}). Note that if $P(\nu,\omega)=0$,
then by Theorem \ref{normality} and the continuous mapping theorem, we get
\begin{equation}
\sqrt{\frac{d_n}{L_{d_n}}}|\hat G_{n}^{c,d}(\nu,\omega)| \stackrel
{d}{\longrightarrow}\mathcal{L}(Z),
\label{z-square}
\end{equation}
where $Z:= \sqrt{S_{1}^2 +S_{2}^2 }$ and the random vector
$(S_{1},S_{2})$ has a two-dimensional normal distribution with
mean zero and a covariance
matrix equal to $\Sigma(\nu,\omega)$. If we assume that $\operatorname{trace}(\Sigma(\nu,\omega))>0,$ then the distribution of
a random variable $Z$ is
continuous. Now, assume that $P(\nu,\omega)\not=0$.
Applying the delta method for the convergence (\ref{eq32})
and the function $f( x,y)=\sqrt{x^2+y^2} $ that is differentiable at
the point $(x_0,y_0)=(\operatorname{Re}(P(\nu,\omega)),\operatorname
{Im}(P(\nu,\omega)))$, we get
\begin{equation} \label{|sigma1|}
\sqrt{\frac{d_n}{L_{d_n}}}\bigl(|\hat G_{n}^{c,d}(\nu,\omega)|-|P(\nu
,\omega)|\bigr)\stackrel{d}{\longrightarrow
}\mathcal{N}_{1}(0,\mathrm{D}_{1}\Sigma(\nu,\omega)\mathrm{D}_{1}^{\mathrm{T}}),
\end{equation}
where
%
\begin{equation}\label{d1}
\mathrm{D}_{1} =\frac{1}{|P(\nu,\omega)|} ( \operatorname{Re}(P(\nu
,\omega)),\operatorname{Im} (P(\nu,\omega)) )
\end{equation}
and $M^\mathrm{T}$ is the transpose of the matrix $M$.

Now, let us establish the asymptotic distribution of a normalized statistic
$|\hat\gamma_{n}^{c,d}(\nu,\omega)|$ for $\nu\not= \omega$ (in the
case $\nu=\omega$, we
have $|\hat\gamma_{n}^{c,d}(\nu,\omega)|=1$). To do
this, we assume that $g_{0}(\nu)g_{0}(\omega) \not=0$.
In the first case, take $\point$ such that $P(\nu,\omega)=0$. By Theorem
\ref{normality}, the continuous mapping theorem,
consistency of the estimators $\hat G_{n}^{c,d}(\nu,\nu)$ and $\hat
G_{n}^{c,d}(\omega,\omega)$ and Slutsky's
lemma, we then have
\begin{eqnarray} \label{11111a}
\sqrt{\frac{d_n}{L_{d_n}}} |\gamma_{n}^{c,d}(\nu,\omega)| &=& \sqrt{\frac
{d_n}{L_{d_n}}}
\frac{ |\hat G_{n}^{c,d}(\nu,\omega)|}{\sqrt{|\Re(\hat
G_{n}^{c,d}(\nu,\nu))\Re( \hat
G_{n}^{c,d}(\omega,\omega))|}}\nonumber\\[-8pt]\\[-8pt]
&\stackrel{d}{\longrightarrow}& \mathcal{L} \biggl(\frac{Z}{\sqrt{g_{0}(\nu
)g_{0}(\omega)}} \biggr).\nonumber
\end{eqnarray}
If the $\operatorname{trace}(\Sigma(\nu,\omega))>0 $, then the limiting
distribution in (\ref{11111a}) is continuous. As a final consideration,
assume that $P(\nu,\omega)\not=0$. In this case, we use the delta
method to show the asymptotic normality of
the normalized estimator $|\hat\gamma_{n}^{c,d}(\nu,\omega)|$. In
doing this, we use convergence (\ref{gamma1}) and the function
$f(x,y,z,t)=\frac{\sqrt{x^2+t^2}}{\sqrt{yz}}$ that is differentiable at
a point
$(x_0,y_0,z_0,t_0)=(\operatorname{Re}(P(\nu,\omega)),g_{0}(\nu
),g_{0}(\omega
),\Im(P(\nu,\omega))).$ Applying the delta
method, we conclude that
\[
\sqrt{\frac{d_n}{L_{d_n}}} \bigl( |\hat\gamma_{n}^{c,d}(\nu,\omega)| -
|\gamma(\nu,\omega)|\bigr)
\stackrel{d}{\longrightarrow}\mathcal{N}_{1} (0, \mathrm{D}_{2} \Psi
(\nu,\omega) \mathrm{D}_{2}^{\mathrm{T}}),
\]
where
%
\begin{equation} \label{Da}
\mathrm{D}_{2}= \frac{|P(\nu,\omega)|}{\sqrt{g_{0}(\nu)g_{0}(\omega)}}
\biggl( \frac{\operatorname{Re}(P(\nu,\omega))}{|P(\nu,\omega)|^2}, -\frac{1}{2
g_{0}(\nu)} ,-\frac{1}{2 g_{0}(\omega)}, \frac{\operatorname{Im}(P(\nu
,\omega))}{|P(\nu,\omega)|^2} \biggr).
\end{equation}

\begin{remark}
By elementary calculation, the assumption that $ \operatorname{trace}(\Sigma(\nu
,\omega))>0 $ is equivalent to
$g_{0}(\nu) g_{0}(\omega) + |P(\nu,2\uppi-\omega)|^2 >0,$ which is true
if we assume that $g_{0}(\xi)>0$ for
any $\xi\in(0,2\uppi]$. Moreover, $\mathrm{D}_{1}\Sigma(\nu,\omega) \mathrm{D}_{1}^{\mathrm{T}}> 0 $ if we assume that
$\operatorname{det}(\Sigma(\nu,\omega))\not= 0$, and $\mathrm{D}_{2} \Psi(\nu
,\omega) \mathrm{D}_{2}^{\mathrm{T}}>0 $ if we
assume that $\operatorname{det}(\tilde\Psi(\nu,\omega)) \not= 0$.
\end{remark}

In the corollaries which follow, we summarize considerations concerning
the asymptotic distributions of normalized estimators $|\hat
G_{n}^{c,d}(\nu,\omega)|$ and $|\hat\gamma_{n}^{c,d}(\nu,\omega)|$.

\begin{corollary}\label{J1-dis}
Let $\ciag{X_{t}}$ be a time series such that all the assumptions of
Theorem \ref{normality}
hold and that $g_{0}(\xi) \not=0$ for any $\xi\in(0,2\uppi]$. Take
any $\point$.
For the case where $P(\nu,\omega)\not=0$, we require that $
\operatorname{det}(\Sigma(\nu,\omega))\not=0$,
where matrix $\Sigma(\nu,\omega)$ is given by Theorem \ref{normality}.
Under these assumptions, we have the convergence
\begin{eqnarray}\label{limitingp}
&&\sqrt{\frac{d_n}{L_{d_n}}}\bigl(|\hat G_{n}^{c,d}(\nu,\omega)| -|P(\nu
,\omega)|\bigr)\nonumber\\[-8pt]\\[-8pt]
&&\quad  \stackrel{d}{\longrightarrow} J^{P(\nu,\omega)} :=
\cases{
\mathcal{L}(Z) &\quad  if $ P(\nu,\omega)=0$,\cr
\mathcal{N}_1 (0,\mathrm{D}_{1}\Sigma(\nu,\omega)\mathrm{D}_{1}^{\mathrm{T}})
&\quad if $ P(\nu,\omega)\not=0$,
}\nonumber
\end{eqnarray}
where $\mathcal{L}(Z)$ and $\mathrm{D}_{1}$ are given by (\ref{z-square})
and (\ref{d1}), respectively. Moreover, the distribution
$J^{P(\nu,\omega)}$ ($J^{P}$ for short) is continuous.
\end{corollary}
%

\begin{corollary}\label{J2-dis-a}
Suppose that all the assumptions of Theorem \ref{normality} hold and
that $g_{0}(\xi) \not=0$ for any $\xi\in(0,2\uppi]$.
Take any $\point$ such that $\nu\not= \omega$. For the case where
$P(\nu,\omega)\not=0$, we
require that $\operatorname{det}( \Psi(\nu,\omega))\not=0$,
where the matrix $ \Psi(\nu,\omega)$ is given by (\ref{gamma1}). Then
\begin{eqnarray}\label{limitingg}
&&\sqrt{\frac{d_n}{L_{d_n}}} \bigl( |\hat\gamma_{n}^{c,d}(\nu,\omega)| -
|\gamma(\nu,\omega)|\bigr)\nonumber\\[-8pt]\\[-8pt]
&&\quad  \stackrel{d}{\longrightarrow}
J^{\gamma(\nu,\omega)} := \cases{
\mathcal{L}\bigl(Z / \sqrt{g_{0}(\nu)g_{0}(\omega)}\bigr) & \quad if $ P(\nu,\omega)=0$,\cr
\mathcal{N}_1 (0,\mathrm{D}_{2}\Psi(\nu,\omega)\mathrm{D}_{2}^{\mathrm{T}}) & \quad if $ P(\nu,\omega)\not=0$,
}\nonumber
\end{eqnarray}
where $\mathcal{L}(Z /\sqrt{g_{0}(\nu)g_{0}(\omega)})$ and $\emph
{D}_{2}$ are given by (\ref{11111a}) and (\ref{Da}), respectively.
Moreover, the distribution $J^{\gamma(\nu,\omega)}$ ($J^{\gamma}$ for
short) is continuous.
\end{corollary}

Note that we cannot use the asymptotic distributions from Corollaries
\ref{J1-dis} and \ref{J2-dis-a} in practice if we
do not know that $P(\nu,\omega)=0$. Even if we know the period $T$ for
PC time series, the value of $|P(\nu,\nu-\lambda)|$
for $\lambda\in \{2 k \uppi/T\dvtx k=0,1,\ldots,T-1\}$ can be equal to
zero. For example, consider the model (\ref{pmaa}) with period
$T=12$.
After calculations (we omit overly long formulae), ${g_0}(\nu) =
(235+72 \cos(\nu))/(16
\uppi)$, ${g_{\uppi/{6}}}(\nu) \not\equiv0$, ${g_{\uppi
/{3}}}(\nu) \not\equiv0$,
${g_{\uppi/{2}}}(\nu) = \frac{1}{2\uppi}$, $ {g_{2 \uppi
/{3}}} (\nu) = \frac{1}{32 \uppi}$
and, crucially for this example, ${g_{5 \uppi/{6}}} (\nu) =
0$, ${g_{\uppi}} (\nu) = 0$ for $\nu\in(0,2\uppi]$. Therefore, we refer
to the subsampling methodology, which does not require information
about $P(\nu,\omega)$.

\section{Subsampling procedure and consistency}\label{sec4}

Let us introduce the subsampling procedure. Some of the notation is
adopted from \cite{Pol}. By $L_{n,b}^{P}$ and
$L_{n,b}^{\gamma}$, we denote the subsampling estimators of
distribution functions of
$\sqrt{n / L_{n}} ( |\hat G_{n}(\nu,\omega)| - |P (\nu,\omega)|)$ and
$ \sqrt{n / L_{n}} ( |\hat\gamma_{n}(\nu,\omega)| - |\gamma(\nu
,\omega)|)$, respectively.
For any point $x \in\mathbb{R}$,
we define those estimators as
\begin{eqnarray}
L_{n,b}^{P}(x)&=&\frac{1}{n-b+1}\sum_{t=1}^{n-b+1}\mathbf{1}\bigl\{\sqrt
{b/L_{b}} \bigl( |\hat G_{n}^{t-1,b}(\nu,\omega)| - |\hat G_{n}(\nu
,\omega)|\bigr)\leq
x\bigr\}, \label{sub1}
\\
L_{n,b}^{\gamma}(x)&=&\frac{1}{n-b+1}\sum_{t=1}^{n-b+1}\mathbf{1}\bigl\{\sqrt
{b/L_{b}} \bigl( |\hat\gamma_{n}^{t-1,b}(\nu,\omega)| - |\hat\gamma
_{n}(\nu
,\omega)|\bigr)\leq x\bigr\}, \label{sub2}
\end{eqnarray}
where
\begin{eqnarray*}
\hat G_{n}^{t-1,b}(\nu,\omega)&=&\frac{1}{2 \uppi b}\sum
_{j_1=t}^{t+b-1}\sum_{j_2=t}^{t+b-1}
H_{L_b}(j_1-j_2)X_{j_{1}}X_{j_2} \mathrm{e}^{-\mathrm{i} \nu j_1}
\mathrm{e}^{\mathrm{i} \omega j_2 },
\\
|\hat\gamma_{n}^{t-1,b}(\nu,\omega)|&=&\frac{|\hat G_{n}^{t-1,b}(\nu
,\omega)|}
{\sqrt{|\Re(\hat G_{n}^{t-1,b}(\nu,\nu))
\Re(\hat G_{n}^{t-1,b}(\omega,\omega
))|}}
\end{eqnarray*}
and $\mathbf{1}{\{B\}}$ is an indicator function of the event $B$.
Denote by $J^{P}(x)$ and $J^{\gamma}(x)$ the distribution functions of
$J^{P}$ and $J^{\gamma}$, respectively.
For any $\alpha\in
(0,1)$ let $c_{n,b}^{P}(1-\alpha)$ and $c_{n,b}^{\gamma}(1-\alpha)$ be
quantiles of the nominal level $1-\alpha$ from subsampling distributions
(\ref{sub1}) and (\ref{sub2}), respectively.
The following theorems concerning consistency of subsampling then hold.

\begin{theorem}\label{subsamplingp}
Under all assumptions of Corollary \ref{J1-dis}:
\begin{longlist}[(iii)]
\item[(i)]$L_{n,b}^{P}(x) \stackrel{p}{\rightarrow} J^{P}(x)$
for any $x \in\mathbb{R}$;
\item[(ii)]$\sup_{x \in
\mathbb{R}}|L_{n,b}^{P}(x) - J^{P}(x)|\stackrel{p}{\rightarrow}0$;
\item[(iii)] the subsampling confidence intervals for
$|P(\nu,\omega)|$ are consistent, which means that
\begin{equation}\label{con-sub1}
P \bigl(\sqrt{n/L_{n}} \bigl(|\hat G_{n}(\nu,\omega)|-|P(\nu,\omega)| \bigr)\leq
c_{n,b}^{P}(1-\alpha) \bigr)\longrightarrow1-\alpha,
\end{equation}
where $b=b(n) \rightarrow\infty$ and $b/n \to0$.
\end{longlist}
\end{theorem}

\begin{theorem}\label{subsamplingc}
Under all the assumptions of Corollary \ref{J2-dis-a}:
\begin{longlist}
\item[(i)]$L_{n,b}^{\gamma}(x) \stackrel{p}{\rightarrow}
J^{\gamma}(x)$ for any $x \in\mathbb{R}$;
\item[(ii)]$\sup_{x \in
\mathbb{R}}|L_{n,b}^{\gamma}(x) - J^{\gamma}(x)|\stackrel{p}{\rightarrow
}0$;
\item[(iii)] the subsampling confidence intervals for
$|\gamma(\nu,\omega)|$ are consistent, which means that
\begin{equation} \label{con-sub1}
P \bigl(\sqrt{n/L_{n}} \bigl(|\hat\gamma_{n}(\nu,\omega)|-|\gamma(\nu,\omega)|
\bigr)\leq c_{n,b}^{\gamma}(1-\alpha) \bigr)\longrightarrow
1-\alpha,
\end{equation}
where $b=b(n) \rightarrow\infty$ and $b/n \to0$.
\end{longlist}
\end{theorem}

\section{Simulation study}\label{sec5}

In this section, we show possible applications of the results from
Sections \ref{sec3} and \ref{sec4}. In Section~\ref{sec5.1}, we
present graphical methods for detecting the presence of periodic
autocorrelation. We consider simulated data sets.
In Section \ref{sec5.2}, we calculate confidence intervals
for the magnitude of a spectral density using subsampling and
asymptotic distributions.

\subsection{Graphical methods for detecting periodicity}\label{sec5.1}

Take any point $(\nu, \omega)$ from the bifrequency square $(0,2\uppi]^2$
such that $\nu\not= \omega$. Let us consider the null hypothesis $H_0
\dvtx |P(\nu, \omega)|
=0$ and the alternative hypothesis $H_1 \dvtx|P(\nu, \omega)| \not=0$. If
we assume that $g_{0}(\nu)g_{0}(\omega) \not= 0$, then the above
hypotheses $H_{0}$ and $H_{1}$ are equivalent to $H_0 \dvtx |\gamma(\nu,
\omega)| =0$ and $H_1 \dvtx|\gamma(\nu, \omega)| \not=0,$ respectively.
Under the assumptions of Theorems \ref{subsamplingp} and \ref{subsamplingc}, we can use the statistics
$\Upsilon_n^{P}(\nu, \omega)=\sqrt{n/L_{n}} |\hat G_{n}(\nu, \omega
)|$ and $\Upsilon_n^{\gamma}(\nu, \omega)=\sqrt{n/L_n} |\hat\gamma_{n}(\nu,
\omega)|,$ and critical values from the subsampling distributions (\ref{sub1}) and
(\ref{sub2}) for the above testing problems. Under $H_{0}$, the rejection
probability tends to $\alpha$ and under $H_{1}$, this probability tends
to one, which means that the both tests
are asymptotically consistent.

Additionally, by Corollary \ref{J1-dis}, we use the statistics
\[
\tilde\Upsilon_n^{P}(\nu, \omega)=nL_{n}^{-1} \bigl(
\bigl(\Re(\hat G_{n}(\nu, \omega))/\hat\sigma_{n}^{R}(\nu,\omega) \bigr)^2 +
\bigl(\Im(\hat G_{n}(\nu, \omega))/\hat\sigma_{n}^{I}(\nu,\omega) \bigr)^2 \bigr)
\]
for the above testing problem,
where $\hat\sigma_{n}^{R}(\nu,\omega)$ and $\hat\sigma_{n}^{I}(\nu
,\omega)$ are estimators of the asymptotic
variances
of $\sqrt{n/L_{n}}\Re(\hat G_{n}(\nu, \omega))$ and $\sqrt{n/L_{n}}\Im
(\hat G_{n}(\nu, \omega))$, respectively
(see Theorem \ref{covariance}),
obtained by replacing the unknown spectral densities by their
estimates (see formula (\ref{dens})). Under hypothesis $H_{0}$, the
matrix $\Sigma(\nu,\omega)$ from Theorem \ref{normality}
has non-zero
elements only on the main diagonal. Therefore, under $H_{0}$, the test
statistic $\tilde\Upsilon_n^{P}(\nu, \omega)$
has,
asymptotically, a chi-square distribution with two degrees of freedom.

For different points $(\nu, \omega)$ from the bifrequency square
$(0,2\uppi]^2$, we may consider the above testing problems. If we reject
hypothesis $H_{0}$ for many points for which $\omega=\nu-
\lambda$, where $\lambda$ is some constant different from zero,
then it can be suspected that the time series $\ciag{X_{t}}$ is not
stationary (see \cite{hurdgra}, Figure 1).

\begin{figure}

\includegraphics{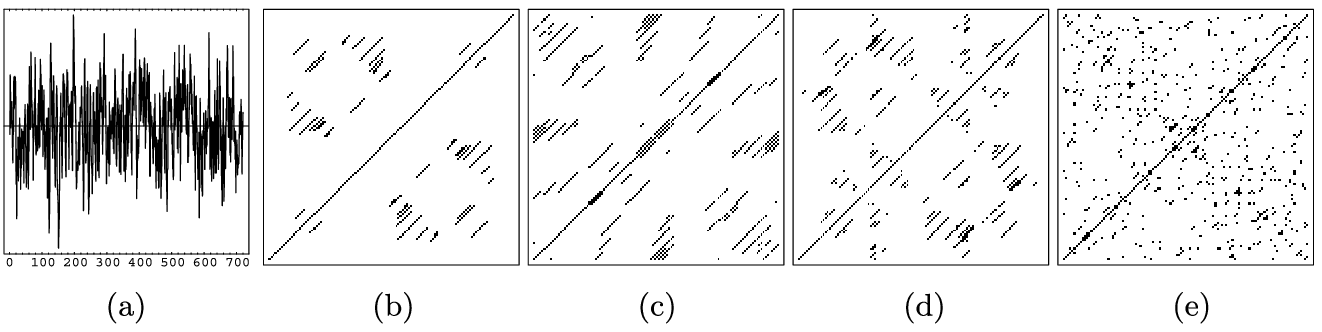}

\caption{Graphical methods for detecting the presence
of periodicity of covariance for samples from the MA(2) model:
$X_{t}=2\epsilon_{t-2} + \epsilon_{t-1} + \epsilon_{t}$, with
$\operatorname{Normal}(0,1)$-distributed innovations: (a) time series;
(b)~tests for $|P(\nu_s,\omega_t)|$ based on the statistic $\Upsilon
_n^{\gamma}(\nu_s, \omega_s)$ and
subsampling; (c) tests for $|P(\nu_s,\omega_t)|$ based on the
statistic $\tilde\Upsilon_n^{P}(\nu_s, \omega_t)$ and
central chi-square distribution with two degrees of freedom; (d) tests
for $|\gamma(\nu_s,\omega_t)|$ based on the statistic
$\Upsilon_n^{\gamma}(\nu_s, \omega_t )$ and subsampling; (e) tests
for $|\gamma(\nu_s,\omega_t)|$ based on the statistic
presented in \protect\cite{hurdgra}, formula (11), page 342, and
the beta distribution.} \label{fig_hist1}
\end{figure}

\begin{figure}[b]

\includegraphics{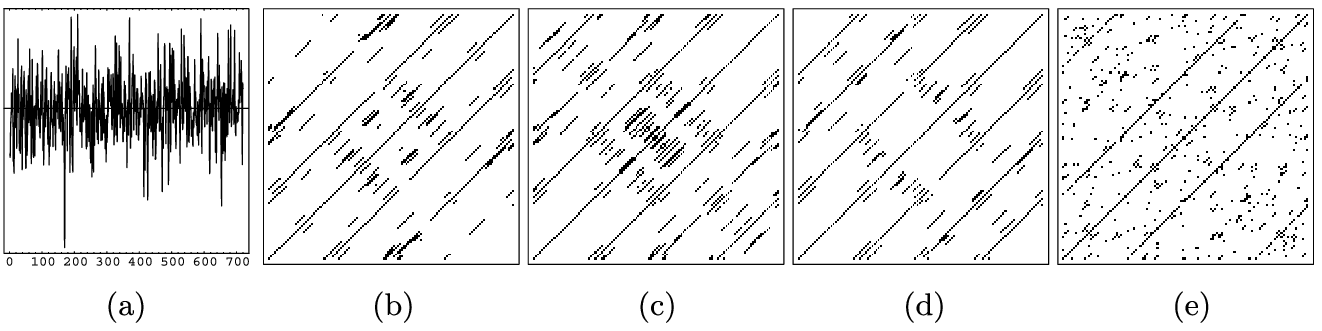}

\caption{Graphical methods for detecting the presence of
periodicity of covariance for samples from the PMA(1) model:
$X_{t}=(2+\sin(2 \uppi t/4))^2\epsilon_{t-1}+\epsilon_{t}$, with period
$T=4$ and
Normal(0,1)-distributed innovations $\epsilon_{t}$: (a) time series; (b)
tests for $|P(\nu_s,\omega_t)|$ based on the statistic $\Upsilon
_n^{\gamma}(\nu_s, \omega_s)$ and
subsampling; (c) tests for $|P(\nu_s,\omega_t)|$ based on the
statistic $\tilde\Upsilon_n^{P}(\nu_s, \omega_t)$ and
central chi-square distribution with two degrees of freedom; (d) tests
for $|\gamma(\nu_s,\omega_t)|$ based on the statistic
$\Upsilon_n^{\gamma}(\nu_s, \omega_t )$ and subsampling; (e) tests
for $|\gamma(\nu_s,\omega_t)|$ based on the statistic
presented in \protect\cite{hurdgra}, formula (11), page 342, and
the beta distribution.} \label{fig_hist2}
\end{figure}

We calculate values of the statistic $\Upsilon_n^{P}(\nu, \omega
)$ for data for each point from the set $\{(\nu_{s},\omega_{t})=(2
\uppi s/120,2 \uppi t/120), 1 \leq s \leq 120, 1 \leq t \leq
120, s \not= t\} \subset(0,2\uppi]^2$ and we compare these values
with the critical values $c_{n,b}^{P}(1-\alpha)$. If the value of
$\Upsilon_n^{P}(\nu, \omega)$ exceeds $c_{n,b}^{P}(1-\alpha)$,
then we reject hypothesis $H_{0}$ and put a black mark at the
point $(\nu,\omega)$ on the bifrequency square. The same steps are
repeated using statistics $\Upsilon_n^{\gamma}(\nu, \omega)$,
$\tilde\Upsilon_n^{P}(\nu, \omega)$ and appropriate critical
values from subsampling and chi-square distributions. We fix
$n=720$, $L_n=[n^{1/5}]=4$, $b=[3 \sqrt{n}]=80$,
$L_{b}=[b^{1/5}]=2$ and $\alpha=99\%$. We put
$H_{L_{n}}(\tau)=\mathbf{1}\{ |\tau| \leq L_{n}\}$, where
$\mathbf{1}{\{B\}}$ is an indicator function of the event $B$. We
compare our methods with the method presented in \cite{hurdgra} and \cite{hurdmiamee}. We put $M=25$ into
\cite{hurdgra}, formula (11). Note that using multiple
testing procedures for
various points seems to be difficult for graphical
methods presented in Figures \ref{fig_hist1}(b)--(e) and \ref{fig_hist2}(b)--(e)
because the number of different tests on each figure is 14\,400 and this
is bigger than the length of the sample, which is $n=720$. If we use
multiple testing
procedures and we calculate the
quantiles from the subsampling distribution,
then the justified significant level
should not be too small in comparison with the sample size. Some
intuitive adjustment for the significant level in the graphical method
presented in Figures \ref{fig_hist1}(e) and  \ref{fig_hist2}(e) for the PC case was made in
\cite{hurdmiamee}. The theoretical results which can
help in
multiple testing procedures in connection with subsampling
approximation can be found in \cite{bertailpolitis01}, where
the question of how well the subsampling distribution approximates
the tail of an unknown distribution that we approximate was
considered.

Figure \ref{fig_hist1} presents the stationary case. All methods
for this case fail to show the evidence of periodic correlation.
Figure \ref{fig_hist2} presents methods for the sample from
the periodic moving average model of order one with period equal
to $T=4$. All methods for this case clearly reveal the presence of
periodic correlation with appropriate length of period $T=4$.
In summary, all methods can successfully distinguish the stationary
case from the PC case. Moreover, there is no clear difference in
these examples between the efficiency of the above graphical methods based
on subsampling, the central chi-square distribution with two
degrees of freedom and the beta distribution.\looseness=-1

We chose the parameters $L_n$ and $b$ to satisfy the assumptions
in Theorems \ref{subsamplingp} and \ref{subsamplingc}. Note that in \cite{Pol}, pages
190--191, the problem of choosing $b$ in practice in special cases was
considered. It was found that the choice of $b$
is, in practice, a delicate issue, like the choice of the
bandwidth $L_n$ in nonparametric spectral density estimation. It
was also shown in \cite{Len08b} in simulations
that the type I
error and power
of the test based on subsampling approximation for significance of the
value $|a(\lambda,\tau)|$
strongly depends on the choice of the
parameter $b$. It is likely that a similar problem occurs in the
graphical methods presented in this work, where the choice of
the parameters $L_{n}$ and $b$ is a open problem.

\subsection{Confidence intervals}\label{sec5.2}

From Theorem \ref{subsamplingp}, we can obtain a two-sided equal-tailed
consistent confidence interval for the parameter
$|P(\nu,\omega)|$ in the form
\begin{equation}\label{sub_c_i}
\biggl(|\hat G_{n}(\nu,\omega)| -c_{n,b}^P\biggl(1-\frac{\alpha}{2}\biggr) \Big/ \sqrt
{n/L_n} , |\hat G_{n}(\nu,\omega)|
-c_{n,b}^P\biggl(\frac{\alpha}{2}\biggr)\Big/ \sqrt{n/L_{n}} \biggr),
\end{equation}
where $\alpha\in(0,1)$ is a nominal level. A similar confidence
interval can be obtained for coherence using Theorem \ref
{subsamplingc}. Consider the periodic moving average
model of order one with period equal to $T$,
\begin{equation} \label{pmaa}
X_{t}=\theta(t-1) \epsilon_{t-1}+\epsilon_{t},
\end{equation}
where $\theta(t)=(2+\sin(2 \uppi t/T))^2$ and $\epsilon_{t}$ is a
Gaussian white noise with mean zero and variance equal to
one. The aim of this section is to calculate confidence intervals for
$|g_{\lambda}(\nu)|$ for time series of the form (\ref{pmaa}).

\begin{figure}
\begin{tabular}{@{}c@{}c@{}}

\includegraphics{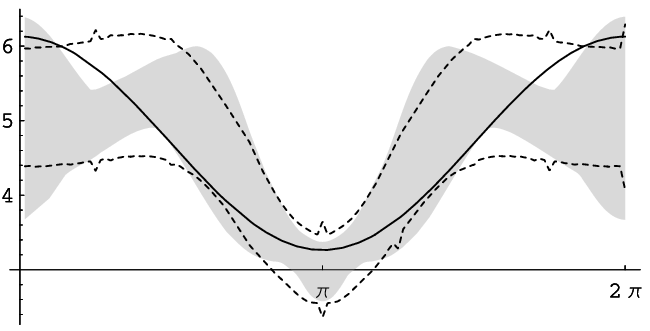}
&\includegraphics{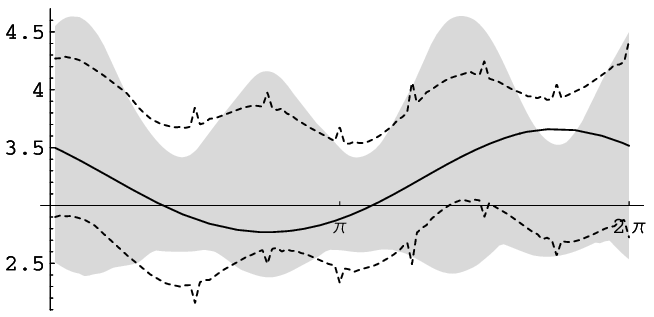}\\
(a)&(b)\\

\includegraphics{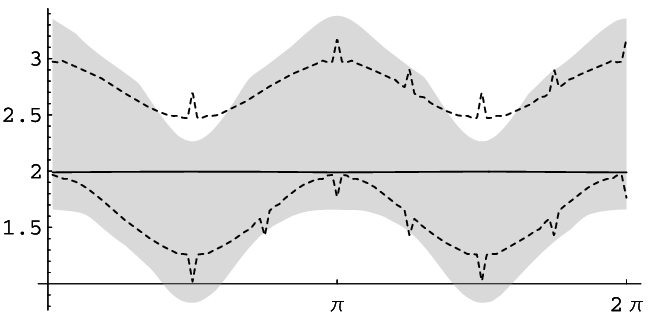}
&\includegraphics{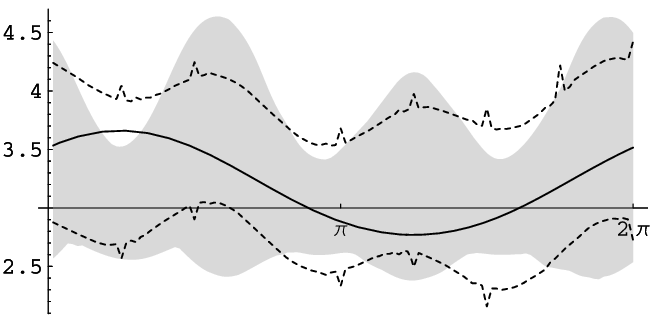}\\
(c)&(d)
\end{tabular}
\caption{Confidence internals for
$|P(\nu,\nu-\lambda)|$ for $\lambda\in\Lambda=\{0,\uppi/2,\uppi,\frac
{3}{2}\uppi\}$ and for $\nu\in\{ 2 k \uppi/120\dvtx k=1,2,\ldots,120\}$.
Solid line: theoretical $|P(\nu,\nu-\lambda)|$;
dashed line: $95\%$ asymptotic confidence interval (based on
distribution $J^{P}$ from Corollary \protect\ref{J2-dis-a})
for sample size $n=500$;
light shading: $95\%$ subsampling confidence interval
(\protect\ref{sub_c_i}). (a) The case $\lambda=0$,
$g_{0}(\nu)=|P(\nu,\nu)|=(59+18 \cos(\nu))/4 \uppi$; (b) the case $\lambda=\frac{\uppi}{2}$,
$|g_{\uppi/2}(\nu)|=|P(\nu,\nu-\uppi/2)|=\sqrt{2} (\sqrt{51+10 \cos
(\nu)-10 \sin(\nu)-\sin(2 \nu)})/\uppi$;
(c) the case $\lambda=\uppi$, $|g_{\uppi}(\nu)|=|P(\nu
,\nu-\uppi)|=\sqrt{625+4 \sin^2(\nu)}/4 \uppi$; (d) the case
$\lambda=\frac{3}{2} \uppi$, $|g_{(3/2)\uppi}(\nu)|=
|P(\nu,\nu-\frac{3}{2} \uppi)|={\sqrt{2}} {\sqrt{51+10 \cos(\nu)+
10 \sin(\nu) + \sin(2 \nu)}} / \uppi$.}\label{fig_hist}
\end{figure}

Figure \ref{fig_hist} presents $95 \%$ confidence intervals (\ref{sub_c_i}) for
$|P(\nu,\nu-\lambda)|$, $T=4$, $\lambda
\in\Lambda=\{0,\uppi/2,\uppi,\frac{3}{2}\uppi\}$ and for $\nu\in\{ 2 k \uppi
/120\dvtx k=1,2,\ldots,120\}$.
We put $L_n=[n^{1/5}]$, $b=[3 \sqrt{n}]$ and $n=500$. We also present
$95\%$ confidence intervals based on the asymptotic distribution $J^{P}$
(see the formula in Corollary \ref{J1-dis}) by replacing unknown values
of spectral densities in the asymptotic distribution
$J^{P}$ by their estimates using (\ref{dens}). The confidence interval
for $g_{\lambda}(\nu)$ based on the asymptotic
distribution $J^{P}$ has spikes
on the interval $(0,2\uppi]$ being considered. It follows from the
fact that asymptotic variance in the limiting distribution $J^P$
depends (for fixed $\lambda$) on a value $P(-\nu,\nu-\lambda)$
that is non-zero for the considered $\nu\in\{ 2 k \uppi/120\dvtx
k=1,2,\ldots,120\}$ only for the
points $\nu\in\{\uppi\pm x/2\dvtx x \in\Lambda\}$. This can also be
observed for the stationary case, where the asymptotic variance of the
normalized spectral density estimator
on the main diagonal, expressed as a function of $\nu$, has a
discontinuity at points $\nu_0=\uppi$ and $\nu_0=2 \uppi$
(see \cite{priestley} and \cite{zurbenko}). Note, that the distribution
$J^P$ can be used to construct a confidence interval for $|P(\nu,\omega
)|$ under additional information concerning parameter $P(\nu,\omega)$.
Therefore, subsampling seems to be more useful in practice when there
is no rule to obtain the parameter $b$ in our testing
problem.
%

%

\section{Conclusions and open problems}\label{sec6}
In this work, we give the exact forms of the asymptotic
distributions for normalized spectral density and magnitude of
coherence statistics and we prove the consistency of subsampling
applied for magnitude of spectral density and coherence. This
research was done for $\alpha$-mixing zero-mean APC time series
with the support of the spectral measure concentrated on a countable set
of parallel lines. By the consistency of subsampling, we construct
asymptotically consistent confidence intervals for magnitude of
spectral density and coherence. Importantly, this
confidence intervals does not depend on the support of the spectral
measure, unlike confidence intervals based on asymptotic
distributions. Graphical methods for detecting the presence of periodic
autocovariance were presented and their theoretical properties
explored. One of the main problems is the choice of the parameter
$b$ in graphical methods, as presented in Section \ref{sec5}.

The next problem connected with applications is the estimation of
spectral density and coherence for non-zero-mean time series.
Note that the well-known differencing operation, popular in
time series analysis, is useful when the mean function is a
constant or a periodic function with known period. Generally, if we
assume that the mean function is almost periodic, then we cannot
use the differencing operation, even if we known the set $\Gamma$ (see
Remark \ref{non-zero}), in the Fourier representation of the mean function.
To cope with this problem, we can consider spectral density
estimators for non-zero-mean time series of the form (\ref{dens2}).
This, and other open problems, are currently being
researched by the author.

\begin{appendix}

\section*{Appendix}\label{app}

\renewcommand{\theequation}{\arabic{equation}}

To begin, we formulate a few auxiliary lemmas which are crucial for
proving the theorems.
\begin{lemma}\label{base}
Let $\ciag{X_{t}}$ be real-valued time series with corresponding $\alpha
$-mixing sequence $\alpha(\cdot)$.
Assume that there exist real numbers $\delta>0$ and $\Delta<\infty$
such that $\sup_{t \in\mathbb{Z}} \| X_t \|_{2+\delta} \leq\Delta$. For any
$\lambda\in[0,2\uppi)$, $j \in\mathbb{Z}$ and
$\tau\in\mathbb{Z}$, we then have the following estimates:
\begin{longlist}[(ii)]
\item[(i)] $|B(j,\tau)| \leq 8 \Delta^2 \alpha^{\delta/{(2
+ \delta)}}(|\tau|)$;
\item[(ii)] $|a(\lambda,\tau)| \leq 8 \Delta^2 \alpha^{
{\delta}/{(2 + \delta)}}(|\tau|)$.
\end{longlist}
\end{lemma}

\begin{pf}
The proof of (i) is elementary and is based on \cite{doukhan}, Theorem
3, page 9. The proof of (ii) follows from (i) by noting
that $a(\lambda,\tau)=\lim_{n \to\infty} n^{-1} \sum_{t=1}^{n}
B(t,\tau) \mathrm{e}^{-\mathrm{i} \lambda t}$ for any
$\lambda\in[0,2\uppi)$ and $\tau\in\mathbb{Z}$.
\end{pf}

\begin{lemma}[(Conclusion from \cite{doukhan}, Theorem 2, page 26)]\label{fazekas}
Let $\ciag{X_{t}}$ be a zero-mean time series. Assume that there exist
real numbers $l>2$ and $\delta>0$ and an even
integer $r \geq l-2$ such that $\sup_{t \in\mathbb{Z}} \norm
{X_{t}}_{l+\delta} < \infty$ and
$\sum_{k=1}^{\infty} k^{r}\alpha^{\delta/{(r+2+\delta)}}(k)<\infty
.$ There then exists a constant $K_{\alpha,l}<\infty$ (which depends
only on the sequence $\alpha(\cdot)$ and constant $l$) such that for
any $a \in
\mathbb{Z}$ and $g \in\mathbb{N}$, we have the estimate
\begin{equation}
\Biggl\| \sum_{t=a+1}^{a+g} X_t \Biggr\|_{l} \leq (K_{\alpha,l} \max\{
Q(l,\delta,a,g),[Q(2,\delta,a,g)]^{l/2} \})^{1/l},
\end{equation}
where $Q(l,\delta,a,g)=\sum_{t=a+1}^{a+g} \| X_{t} \|_{l+\delta}^{l}$.
\end{lemma}

\begin{lemma}[(Multidimensional generalization of the CLT of \cite{francq})]\label{as-norm}
Let $\{Y_{n,t}=(y_{n,t,1},y_{n,t,2},\ldots,y_{n,t,d})\in\mathbb{R}^d\dvtx
1 \leq t \leq a_n\}$ be a triangular array
of zero-mean real-valued random vectors, where $\{a_n \}_{n \in\mathbb
{N}}$ is a sequence of positive
integers tending to infinity. Define
\[
S_{n}=\frac{1}{\sqrt{a_n}}\sum_{j=1}^{a_n}Y_{n,j}.
\]
Assume that:
\begin{enumerate}[(iii)]
\item[(i)] there exists a constant $\delta> 0$ such that
\[
\sup\{ \| y_{n,t,k} \|_{2+\delta} \dvtx n \geq1, 1 \leq t \leq a_n, 1 \leq
k \leq d \} <\Delta< \infty;
\]
\item[(ii)] $\lim_{n \to\infty} \operatorname{Var}(S_{n})=\Sigma;$
\item[(iii)] there exists a sequence of positive integers $\{
h_n\}_{n \in\mathbb{N}}$ such that $h_n=\mathrm{O}(a_n^{\kappa})$ for some
$\kappa\in[0,\frac{\delta}{4+4 \delta})$, and a sequence $\{
\varsigma(k) \}_{k\geq0}$ such that $\alpha_{n}(k) \leq
\varsigma(k-h_n)$ for all $k \geq h_n$ and
\[
\sum_{k=1}^{\infty}k^{r}\varsigma(k)^{\delta/{(2+\delta)}}<K<\infty
\]
for some $r> \frac{2 \kappa(1+\delta)}{\delta-2 \kappa(1+\delta)}$,
where $\alpha_{n}(\cdot)$ is a mixing sequence for a triangular array
$\{Y_{n,t} \dvtx 1 \leq t \leq a_n\}$ defined for any $n\geq2$ and
$h=1,2,\ldots,a_n-1$ as
\[
\alpha_{n}(h)= \sup_{1 \leq t \leq a_n-h }\mathop{ \sup_{A \in A_{n,t}}}_{ B \in B_{n,t+h}} |P(A \cap B)-P(A)P(B)|.
\]
$A_{n,t}$ stands for the $\sigma$-algebra generated by $\{Y_{n,u}\dvtx 1
\leq u \leq t\}$ and $B_{n,t}$ stands for the $\sigma$-algebra
generated by $\{Y_{n,u}\dvtx t \leq u \leq a_n\}$.
\end{enumerate}
Then $ S_{n} \stackrel{d}{\longrightarrow} \mathcal{N}_{d}(0,\Sigma).$
\end{lemma}

\begin{pf}
The proof is based on the same steps as that of the CLT of \cite{francq}, and the Cram\'er--Wold device.
\end{pf}

\begin{lemma}\label{inequality1}
Let $\ciag{X_{t}}$ be a zero-mean and $\alpha$-mixing time series.
Assume that there exists a real number $\delta>0$ such that:
\begin{longlist}[(ii)]
\item[(i)] $\sup_{t \in\mathbb{Z}} \norm{X_{t}}_{6+3\delta
}\leq\Delta<\infty;$
\item[(ii)] $\sum_{k=0}^{\infty}
(k+1)^2\alpha(k)^{\delta/{(2+\delta)}}\leq K<\infty$.
\end{longlist}
For any sequences $\{c_n\}_{n \in\mathbb{Z}}$ and $\{d_n\}_{n \in
\mathbb{Z}}$ of integers where $d_n>0$, we then have
\begin{eqnarray*}
&&\sum_{s=c_n +1}^{c_n+d_n} \sum_{t=c_n+1}^{c_n +d_n} \sum_{u=c_n+1}^{c_n+d_n}
\sum_{v=c_n+1}^{c_n+d_n}  |
E(X_{s}X_{t}X_{u}X_{v})-E(X_{s}X_{t})E(X_{u}X_{v}) \\
&&\hphantom{\sum_{s=c_n +1}^{c_n+d_n} \sum_{t=c_n+1}^{c_n +d_n} \sum_{u=c_n+1}^{c_n+d_n}
\sum_{v=c_n+1}^{c_n+d_n}  |}
{} - E(X_{s}X_{u})E(X_{t}X_{v}) - E(X_{s}X_{v})E(X_{t}X_{u}) | \leq C
d_n ,
\end{eqnarray*}
where the constant $C$ depends only on $K$ and $\Delta$.
\end{lemma}

\begin{pf} This proof is based on a similar technique as in the
proof of \cite{Len08a}, Theorem 4.1, and is
therefore omitted.
\end{pf}

\begin{lemma}\label{bias}
Suppose that assumptions \textup{(A1)--(A5)} hold. Using the
notation from Section \ref{sec3}, for any $(\nu,\omega) \in(0,2\uppi]^2$, we have
\[
\lim_{n \to\infty} E [\hat P_{n}^{c,d}(\nu,\omega)] = P(\nu,\omega).
\]
%
\end{lemma}

\begin{pf}
Take any $\point$ and note that ($j=t$, $\tau=s-t$)
\begin{eqnarray*}
E[ \hat P_{n}^{c,d}(\nu,\omega)]
&=& \frac{1}{2\uppi d_n}\sum_{s=c_n+1}^{c_n+d_n}\sum
_{t=c_n+1}^{c_n+d_n}E(X_{s}X_{t})\mathrm{e}^{-\mathrm{i} (\nu-\omega)
t}\mathrm{e}^{-\mathrm{i} \nu (s-t)}\\
& =& \frac{1}{2\uppi d_n} \sum_{j=c_n+1}^{c_n+d_n}\sum_{\tau=c_n+1-j}^{c_n+d_n-j}
B(j,\tau)\mathrm{e}^{-\mathrm{i} (\nu-\omega) j}\mathrm{e}^{-\mathrm{i}
\nu \tau}\\
& = &\frac{1}{2\uppi d_n} \sum_{j=c_n+1}^{c_n+d_n}
\sum_{\tau=c_n+1-j}^{c_n+d_n-j}\sum_{\lambda\in\Lambda_\tau
}a(\lambda,\tau)
\mathrm{e}^{\mathrm{i} (\lambda-(\nu-\omega)) j} \mathrm{e}^{-\mathrm
{i} \nu \tau}
\\
&=& \frac{1}{2\uppi d_n} \sum_{j=c_n+1}^{c_n+d_n}\sum_{\tau
=c_n+1-j}^{c_n+d_n-j} a(\nu-\omega,\tau)
\mathrm{e}^{-\mathrm{i} \nu \tau}\\
&&{}+ \frac{1}{2\uppi d_n} \sum_{j=c_n+1}^{c_n+d_n}\sum_{\tau=c_n+1-j}^{c_n+d_n-j}
\sum_{\lambda\in\Lambda_\tau\setminus\{\nu-\omega\} }a(\lambda
,\tau)\mathrm{e}^{\mathrm{i} (\lambda-(\nu-\omega)) j}
\mathrm{e}^{-\mathrm{i} \nu \tau}
\\
&=& \frac{1}{2\uppi} \sum_{|\tau|<d_n} \biggl( 1- \frac{|\tau|}{d_n} \biggr)a(\nu
-\omega,\tau) \mathrm{e}^{-\mathrm{i} \nu \tau}\\
&&{}
+ \frac{1}{2\uppi d_n} \sum_{j=c_n+1}^{c_n+d_n}\sum_{\tau=c_n+1-j}^{c_n+d_n-j}
\sum_{\lambda\in\Lambda_\tau\setminus\{\nu-\omega\} }a(\lambda
,\tau)\mathrm{e}^{\mathrm{i} (\lambda-(\nu-\omega)) j}
\mathrm{e}^{-\mathrm{i} \nu \tau}.
\end{eqnarray*}
Denote the first and second terms of the last equality by $\epsilon
_{1,n}$ and $\epsilon_{2,n}$, respectively.
Note that $\epsilon_{1,n}$ is a Ces\'aro mean
for the sequence $\sum_{|\tau|<d_n} a(\nu-\omega,\tau) \mathrm
{e}^{-\mathrm{i} \nu \tau} $. Hence, $\epsilon_{1,n}$
goes to $P(\nu,\omega)$.
By the estimate $|\sum_{j=p}^{q} \mathrm{e}^{-\mathrm{i} x j}| \leq
|\cosec(x/2)|$, where $p\leq q$,
$x \not\equiv0 \mbox{ modulo } 2 \uppi$, we have
%
\begin{eqnarray*}
|\epsilon_{2,n}| & =& \frac{1}{2\uppi d_n} \Biggl| \Biggl(\sum_{\tau=-d_n+1}^{0} \sum
_{j=c_n+1-\tau}^{c_n+d_n}
+ \sum_{\tau=1}^{d_n-1} \sum_{j=c_n+1}^{c_n+d_n-\tau} \Biggr) \sum
_{\lambda\in\Lambda_\tau
\setminus\{\nu-\omega\} }a(\lambda,\tau)\mathrm{e}^{\mathrm{i} (\lambda
-(\nu-\omega)) j}
\mathrm{e}^{-\mathrm{i} \nu \tau} \Biggr|\\
& \leq& \frac{1}{2\uppi d_n} \sum_{\tau=-d_n+1}^{d_n-1}
\sum_{\lambda\in\Lambda_\tau\setminus\{\nu-\omega\} }
\bigl|a(\lambda,\tau) \cosec\bigl(\bigl(\lambda- (\nu-\omega)\bigr)/2\bigr)\bigr|
\leq \frac{1}{2
\uppi d_n} \sum_{\tau=-d_n+1}^{d_n-1} z_{\tau}(\nu-\omega),
\end{eqnarray*}
which means that $\epsilon_{2,n} \to0$ (by (A3)). This completes the proof.
\end{pf}

\begin{lemma} \label{biasx}
Suppose that \textup{(A1)--(A5)} and \textup{(B)} hold. Additionally,
assume that there
exist $\delta>0$, $\Delta< \infty$ and $K < \infty$ such that:
\begin{longlist}[(ii)]
\item[(i)]
$\sup_{t \in\mathbb{Z}} \norm{X_{t}}_{6+3\delta} <
\Delta;$
\item[(ii)]
$\sum_{k=0}^{\infty}(k+1)^2\alpha(k)^{\delta/{(2+\delta)}}\leq K.$
\end{longlist}
Using the notation from Section \ref{sec3}, for any $(\nu,\omega) \in(0,2\uppi
]^2$ and any $a \in\mathbb{Z}$, we then have
\begin{eqnarray*}
&&\frac{1}{2 \uppi L_{d_n}} \sum_{\tau_1=-L_{d_n}}^{L_{d_n}}
\sum_{\tau_2=-L_{d_n}}^{L_{d_n}} E(X_{\tau_1+a} X_{\tau_2+a})
H_{{L_{d_n}}}(\tau_1)^2 \mathrm{e}^{-\mathrm{i} \nu(\tau_1+a)}
\mathrm{e}^{\mathrm{i} \omega(\tau_2+a)}\\
&&\quad  = \int_{-1}^{1}w(x)^2 \,\mathrm{d}\, x P(\nu,\omega) + e_n(\nu,\omega),
\end{eqnarray*}
where $|e_n(\nu,\omega)|\leq\tilde e_n (\nu,\omega) \to0$ and $\tilde
e_n(\nu,\omega)$ does not depend on $a$.
\end{lemma}

\begin{pf}
Take any $\point$ and note that using the same steps as in the proof of
Lemma~\ref{bias}, we get\vspace*{-1pt}
\begin{eqnarray*}
&& \frac{1}{2\uppi L_{d_n}} \sum_{\tau_1=-L_{d_n}}^{L_{d_n}}
\sum_{\tau_2=-L_{d_n}}^{L_{d_n}} E(X_{\tau_1+a} X_{\tau_2+a})
H_{{L_{d_n}}}(\tau_1)^2 \mathrm{e}^{-\mathrm{i} \nu(\tau_1+a)}
\mathrm{e}^{\mathrm{i} \omega(\tau_2+a)}\\
&&\quad  = \frac{1}{2\uppi L_{d_n}} \sum_{j=a-L_{d_n}}^{a+L_{d_n}}\sum_{\tau
=a-L_{d_n}-j}^{a+L_{d_n}-j}
H_{{L_{d_n}}}^2(j-a)
a(\nu-\omega,\tau) \mathrm{e}^{-\mathrm{i} \nu \tau}\\
&&\qquad {} + \frac{1}{2\uppi L_{d_n}} \sum_{j=a-L_{d_n}}^{a+L_{d_n}}\sum_{\tau
=a-L_{d_n}-j}^{a+L_{d_n}-j}
\sum_{\lambda\in\Lambda_\tau\setminus\{\nu-\omega\} }
H_{{L_{d_n}}}^2(j-a)
a(\lambda,\tau)\mathrm{e}^{\mathrm{i} (\lambda-(\nu-\omega)) j} \mathrm
{e}^{-\mathrm{i} \nu \tau}\\
&&\quad  = \sum_{|\tau|<2 L_{d_n}} \sum_{j=-L_{d_n}}^{L_{d_n}}
\frac{w^2 (j/{L_{d_n}} )}{2\uppi L_{d_n}} a(\nu-\omega,\tau)
\mathrm{e}^{-\mathrm{i} \nu \tau}- \sum_{|\tau|<2 L_{d_n}}
\sum_{j=L_{d_n}-|\tau|}^{L_{d_n}}
\frac{w^2 (j/{L_{d_n}} )}{2\uppi L_{d_n}}\\
&&\hphantom{\quad  = \sum_{|\tau|<2 L_{d_n}} \sum_{j=-L_{d_n}}^{L_{d_n}}
\frac{w^2 (j/{L_{d_n}} )}{2\uppi L_{d_n}} a(\nu-\omega,\tau)
\mathrm{e}^{-\mathrm{i} \nu \tau}- \sum_{|\tau|<2 L_{d_n}}
\sum_{j=L_{d_n}-|\tau|}^{L_{d_n}}}
{}\times a(\nu-\omega,\tau)
\mathrm{e}^{-\mathrm{i} \nu \tau}\\
&&\qquad {} + \frac{1}{2\uppi L_{d_n}} \sum_{\tau=1}^{2 L_{d_n}} \sum
_{j=a-L_{d_n}}^{a+L_{d_n}-\tau}
\sum_{\lambda\in\Lambda_\tau\setminus\{\nu-\omega\} }
w^2 \biggl(\frac{j-a}{L_{d_n}} \biggr) a(\lambda,\tau)\mathrm{e}^{\mathrm{i}
(\lambda-(\nu-\omega)) j} \mathrm{e}^{-\mathrm{i} \nu \tau}\\
&&\qquad {} + \frac{1}{2\uppi L_{d_n}}
\sum_{\tau=-2 L_{d_n}}^{0} \sum_{j=a-L_{d_n}-\tau}^{a+L_{d_n}}
\sum_{\lambda\in\Lambda_\tau\setminus\{\nu-\omega\} }
w^2 \biggl(\frac{j-a}{L_{d_n}} \biggr) a(\lambda,\tau)\mathrm{e}^{\mathrm{i}
(\lambda-(\nu-\omega)) j} \mathrm{e}^{-\mathrm{i} \nu \tau}\\
&&\quad  = u_{0}(n) - u_{1}(n) + u_{2}(n) + u_{3}(n),
\end{eqnarray*}
where $u_{0}(n)$, $u_{1}(n)$, $u_{2}(n)$ and $u_{3}(n)$ are the first,
second, third and fourth terms, respectively, of the
last equation. It it easy to see that $u_{0}(n) \to \int
_{-1}^{1}w^2(x)\,\mathrm{d}x\, P(\nu,\omega)$. Therefore, to
prove the theorem, it is sufficient to show that $u_{1}(n) \to0$,
$u_{2}(n)\to0$ and $u_{3}(n) \to0$. By Lemma
\ref{base}(ii), for the term $u_{1}(n)$, we have
\[
|u_{1}(n)| \leq \sum_{|\tau|<2 L_{d_n}}
\frac{|\tau|+1}{2\uppi L_{d_n}} 8 \Delta^2 \alpha^{\delta/{(2+\delta)
}}(|\tau|)
\leq\frac{32 \Delta^2 K}{2\uppi L_{d_n}} \to0.
\]
For the term $u_{2}(n)$, we have
\begin{eqnarray*}
|u_{2}(n)| & = &\Biggl| \frac{1}{2\uppi L_{d_n}} \sum_{\tau=1}^{2 L_{d_n}}
\sum_{j=-L_{d_n}}^{L_{d_n}-\tau}
\sum_{\lambda\in\Lambda_\tau\setminus\{\nu-\omega\} }
w^2 \biggl(\frac{j}{L_{d_n}} \biggr) a(\lambda,\tau)\mathrm{e}^{\mathrm{i} (\lambda
-(\nu-\omega)) (j+a)} \mathrm{e}^{-\mathrm{i} \nu \tau} \Biggr|\\
& \leq &\Biggl| \frac{1}{2\uppi L_{d_n}} \sum_{\tau=1}^{L_{d_n}}
\sum_{\lambda\in\Lambda_\tau\setminus\{\nu-\omega\} }
\Biggl(\sum_{j=-L_{d_n}}^{-1} + \sum_{j=0}^{L_{d_n}-\tau} \Biggr)
w^2 \biggl(\frac{j}{L_{d_n}} \biggr) \mathrm{e}^{\mathrm{i} (\lambda-(\nu-\omega))
(j+a)} a(\lambda,\tau) \mathrm{e}^{-\mathrm{i} \nu \tau} \Biggr|\\
&&{} + \Biggl| \frac{1}{2\uppi L_{d_n}} \sum_{\tau=L_{d_n}+1}^{2L_{d_n}}
\sum_{\lambda\in\Lambda_\tau\setminus\{\nu-\omega\} }
\sum_{j=-L_{d_n}}^{L_{d_n}-\tau}
w^2 \biggl(\frac{j}{L_{d_n}} \biggr) \mathrm{e}^{\mathrm{i} (\lambda-(\nu-\omega
))(j+a)} a(\lambda,\tau) \mathrm{e}^{-\mathrm{i} \nu \tau} \Biggr|\\
& \leq& \frac{1}{2\uppi L_{d_n}} \sum_{\tau=1}^{L_{d_n}}
\sum_{\lambda\in\Lambda_\tau\setminus\{\nu-\omega\} }
\underbrace{ \Biggl| \sum_{j=-L_{d_n}}^{-1}
w^2 \biggl(\frac{j}{L_{d_n}} \biggr) \mathrm{e}^{\mathrm{i} (\lambda-(\nu-\omega))
j} \Biggr| }_{f_{1}(n)} |a(\lambda,\tau)|\\
&&{} + \frac{1}{2\uppi L_{d_n}} \sum_{\tau=1}^{L_{d_n}}
\sum_{\lambda\in\Lambda_\tau\setminus\{\nu-\omega\} }
\underbrace{ \Biggl| \sum_{j=0}^{L_{d_n}-\tau}
w^2 \biggl(\frac{j}{L_{d_n}} \biggr) \mathrm{e}^{\mathrm{i} (\lambda-(\nu-\omega))
j} \Biggr| }_{f_{2}(n)} |a(\lambda,\tau)|\\
&&{} + \frac{1}{2\uppi L_{d_n}} \sum_{\tau=L_{d_n}+1}^{2L_{d_n}}
\sum_{\lambda\in\Lambda_\tau\setminus\{\nu-\omega\} }
\underbrace{ \Biggl| \sum_{j=-L_{d_n}}^{L_{d_n}-\tau}
w^2 \biggl(\frac{j}{L_{d_n}} \biggr) \mathrm{e}^{\mathrm{i} (\lambda-(\nu-\omega))
j} \Biggr| }_{f_{3}(n)} |a(\lambda,\tau)|.
\end{eqnarray*}
Now, using the property $|\sum_{j=p}^{q} c_j \mathrm{e}^{\mathrm{i} j
x}| \leq c_p |\cosec(x/2)|,$ which is true for any
$x \not\equiv0 $ modulo $ 2\uppi$ and $c_p \geq c_{p+1} \geq
\cdots\geq c_q$
(see \cite{edwards}, Exercise 1.2, page 10) for the terms $f_1(n)$,
$f_2(n)$, $f_3(n)$, we obtain the estimate
\begin{eqnarray*}
|u_{2}(n)| & \leq & \frac{1}{ \uppi L_{d_n}} \sum_{\tau=1}^{L_{d_n}}
\sum_{\lambda\in\Lambda_\tau\setminus\{\nu-\omega\} }
|a(\lambda,\tau)| \bigl|\cosec\bigl(\bigl(\lambda-(\nu-\omega)\bigr)/2\bigr)\bigr|\\
%
%
&&{} + \frac{1}{2\uppi L_{d_n}} \sum_{\tau=L_{d_n}+1}^{2L_{d_n}}
\sum_{\lambda\in\Lambda_\tau\setminus\{\nu-\omega\} }
|a(\lambda,\tau)| \bigl|\cosec\bigl(\bigl(\lambda-(\nu-\omega)\bigr)/2\bigr)\bigr|
\\
&\leq&\frac{1}{ \uppi L_{d_n}} \sum_{\tau=1}^{2 L_{d_n}} z_{\tau}(\nu
-\omega) \to0.
\end{eqnarray*}
%
Using the same arguments, we get that $u_{3}(n) \to0$. This completes
the proof.
\end{pf}

\begin{pf*}{Proof of Theorem \ref{covariance}}
Using simple decomposition, we have
\begin{eqnarray*}
&& \frac{(2\uppi)^2 d_n}{L_{d_n}} \cov (\hat G_{n}^{c,d}(\nu_1 ,\omega_1
),\hat G_{n}^{c,d}(\nu_2 ,\omega_2 )
)\\
&&\quad  = \frac{1}{d_n L_{d_n}} \sum_{s=c_n +1}^{c_n+d_n} \sum
_{t=c_n+1}^{c_n +d_n}
\sum_{u=c_n+1}^{c_n+d_n} \sum_{v=c_n+1}^{c_n+d_n}
\cov(X_s X_t, X_u X_v)H_{L_{d_n}}(s-t)H_{L_{d_n}}(u-v)\\
&&\hphantom{\quad  = \frac{1}{d_n L_{d_n}} \sum_{s=c_n +1}^{c_n+d_n} \sum
_{t=c_n+1}^{c_n +d_n}
\sum_{u=c_n+1}^{c_n+d_n} \sum_{v=c_n+1}^{c_n+d_n}}
{}\times \mathrm{e}^{-\mathrm{i}(\nu_1 s - \omega_1 t - \nu_2 u + \omega_2 v)}\\
&&\quad  = \frac{1}{d_n L_{d_n}} \sum_{s=c_n +1}^{c_n+d_n} \sum
_{t=c_n+1}^{c_n +d_n}
\sum_{u=c_n+1}^{c_n+d_n} \sum_{v=c_n+1}^{c_n+d_n}
\bigl( \cov(X_s X_t, X_u X_v) - E(X_s X_u)E(X_t X_v)\\
&&\hphantom{\quad  = \frac{1}{d_n L_{d_n}} \sum_{s=c_n +1}^{c_n+d_n} \sum
_{t=c_n+1}^{c_n +d_n}
\sum_{u=c_n+1}^{c_n+d_n} \sum_{v=c_n+1}^{c_n+d_n}
\bigl(}
{} - E(X_s X_v)E(X_t X_u)
\bigr)\\
&&\hphantom{\quad  = \frac{1}{d_n L_{d_n}} \sum_{s=c_n +1}^{c_n+d_n} \sum
_{t=c_n+1}^{c_n +d_n}
\sum_{u=c_n+1}^{c_n+d_n} \sum_{v=c_n+1}^{c_n+d_n}}
{}\times H_{L_{d_n}}(s-t)H_{L_{d_n}}(u-v)\mathrm{e}^{-\mathrm{i}(\nu_1 s -
\omega_1 t - \nu_2 u + \omega_2 v)}
\\
&&
\qquad {} + \frac{1}{d_n L_{d_n}}\sum_{s=c_n +1}^{c_n+d_n} \sum
_{t=c_n+1}^{c_n +d_n}
\sum_{u=c_n+1}^{c_n+d_n} \sum_{v=c_n+1}^{c_n+d_n} E(X_s X_u)E(X_t
X_v)H_{L_{d_n}}(s-t)\\
&&\hphantom{\qquad {} + \frac{1}{d_n L_{d_n}}\sum_{s=c_n +1}^{c_n+d_n} \sum
_{t=c_n+1}^{c_n +d_n}
\sum_{u=c_n+1}^{c_n+d_n} \sum_{v=c_n+1}^{c_n+d_n}}
{}\times H_{L_{d_n}}(u-v)\mathrm{e}^{-\mathrm{i}(\nu_1 s -
\omega_1 t - \nu_2 u + \omega_2 v)}
\\
&&
\qquad {} + \frac{1}{d_n L_{d_n}} \sum_{s=c_n +1}^{c_n+d_n} \sum_{t=c_n+1}^{c_n +d_n}
\sum_{u=c_n+1}^{c_n+d_n} \sum_{v=c_n+1}^{c_n+d_n} E(X_s X_v)E(X_t
X_u)\\
&&\hphantom{\qquad {} + \frac{1}{d_n L_{d_n}} \sum_{s=c_n +1}^{c_n+d_n} \sum_{t=c_n+1}^{c_n +d_n}
\sum_{u=c_n+1}^{c_n+d_n} \sum_{v=c_n+1}^{c_n+d_n}}
{}\times H_{L_{d_n}}(s-t)H_{L_{d_n}}(u-v)\mathrm{e}^{-\mathrm{i}(\nu_1 s -
\omega_1 t - \nu_2 u + \omega_2 v)} \\
&&\quad  = s_{1}(n) + s_{2}(n) + s_{3}(n),
\end{eqnarray*}
where $s_{1}(n)$, $s_{2}(n)$, $s_{3}(n)$ are the first, second and the
third term, respectively, of the last equation.
To prove the theorem, it is sufficient to show that $s_{1}(n)\to0$,
$s_{2}(n) \to(2\uppi)^2 P(\nu_1,\nu_2)\overline{P(\omega_1,\omega_2)}$,
$s_{3}(n) \to(2\uppi)^2 P(\nu_1,2\uppi- \omega_1)
\overline{P(\nu_2,2\uppi- \omega_2)}$ as $n \to\infty.$ Note that by
the estimate
\begin{eqnarray*}
|s_{1}(n)| & \leq &\frac{1}{d_n L_{d_n}} \sum_{s=c_n +1}^{c_n+d_n}
\sum_{t=c_n+1}^{c_n +d_n} \sum_{u=c_n+1}^{c_n+d_n}
\sum_{v=c_n+1}^{c_n+d_n} | \cov(X_s X_t, X_u X_v)-E(X_s X_u)E(X_t
X_v)\\
&&\hphantom{\frac{1}{d_n L_{d_n}} \sum_{s=c_n +1}^{c_n+d_n}
\sum_{t=c_n+1}^{c_n +d_n} \sum_{u=c_n+1}^{c_n+d_n}
\sum_{v=c_n+1}^{c_n+d_n} |}
{}-E(X_s X_v)E(X_t X_u) |
\end{eqnarray*}
and Lemma \ref{inequality1}, we get that $s_{1}(n) \to0$. Considering
the term $s_{2}(n)$, we get ($t=j_1$, $s=j_1+\tau_1$, $v=j_2$, $u=j_2+\tau_2$)
\begin{eqnarray*}
s_{2}(n)
& =& \frac{1}{d_n L_{d_n}}\sum_{s=c_n +1}^{c_n+d_n} \sum
_{t=c_n+1}^{c_n +d_n}
\sum_{u=c_n+1}^{c_n+d_n}
\sum_{v=c_n+1}^{c_n+d_n} E(X_s X_u)E(X_t X_v)\\
&&\hphantom{\frac{1}{d_n L_{d_n}}\sum_{s=c_n +1}^{c_n+d_n} \sum
_{t=c_n+1}^{c_n +d_n}
\sum_{u=c_n+1}^{c_n+d_n}
\sum_{v=c_n+1}^{c_n+d_n}}
{}\times H_{L_{d_n}}(s-t)H_{L_{d_n}}(u-v)
\mathrm{e}^{-\mathrm{i}(\nu_1 s - \omega_1 t - \nu_2 u + \omega_2 v)}\\
& = &\frac{1}{d_n L_{d_n}}
\Biggl(
\sum_{j_1=c_n+1}^{c_n+d_n} \sum_{\tau_1=-L_{d_n}}^{L_{d_n}}
- \sum_{\tau_1=-L_{n}}^{-1}\sum_{j_1=c_n+1}^{c_n-\tau_1}
- \sum_{\tau_1=1}^{L_{n}}\sum_{j_1=c_n+d_n-\tau_1+1}^{c_n+d_n}
\Biggr)\\
&&{}\times \Biggl(
\sum_{j_2=c_n+1}^{c_n+d_n} \sum_{\tau_2=-L_{d_n}}^{L_{d_n}}
- \sum_{\tau_2=-L_{n}}^{-1}\sum_{j_2=c_n+1}^{c_n-\tau_2}
- \sum_{\tau_2=1}^{L_{n}}\sum_{j_2=c_n+d_n-\tau_2+1}^{c_n+d_n}
\Biggr)\\
&&{}\times E(X_{j_1+\tau_1}X_{j_2+\tau_2})E(X_{j_1}X_{j_2}) H_{L_{d_n}}(\tau
_1)H_{L_{d_n}}(\tau_2)
\mathrm{e}^{-\mathrm{i}(\nu_1 (j_1+\tau_1) - \omega_1 j_1 - \nu_2
(j_2+\tau_2) + \omega_2 j_2)}\\
& =& d_{0}(n)- d_{1}(n) - d_{2}(n)+d_{3}(n),
\end{eqnarray*}
where
\begin{eqnarray*}
d_0(n) & =& \frac{1}{d_n L_{d_n}} \sum_{j_1=c_n+1}^{c_n+d_n} \sum_{\tau
_1=-L_{d_n}}^{L_{d_n}}
\sum_{j_2=c_n+1}^{c_n+d_n} \sum_{\tau_2=-L_{d_n}}^{L_{d_n}} \varphi
_{j_1 j_2 \tau_1 \tau_2},\\
d_1(n) & =& \frac{1}{d_n L_{d_n}}
\sum_{j_1=c_n+1}^{c_n+d_n} \sum_{\tau_1=-L_{d_n}}^{L_{d_n}}
\sum_{\tau_2=-L_{n}}^{-1}\sum_{j_2=c_n+1}^{c_n-\tau_2} \varphi_{j_1
j_2 \tau_1 \tau_2}\\
&&{}+ \frac{1}{d_n L_{d_n}} \sum_{j_1=c_n+1}^{c_n+d_n} \sum_{\tau
_1=-L_{d_n}}^{L_{d_n}}
\sum_{\tau_2=1}^{L_{n}}\sum_{j_2=c_n+d_n-\tau_2+1}^{c_n+d_n} \varphi
_{j_1 j_2 \tau_1 \tau_2}, \\
d_2(n) & =& \frac{1}{d_n L_{d_n}}
\sum_{j_2=c_n+1}^{c_n+d_n} \sum_{\tau_2=-L_{d_n}}^{L_{d_n}}
\sum_{\tau_1=-L_{n}}^{-1}\sum_{j_1=c_n+1}^{c_n-\tau_2} \varphi_{j_1
j_2 \tau_1 \tau_2}\\
&&{}
+ \frac{1}{d_n L_{d_n}} \sum_{j_2=c_n+1}^{c_n+d_n} \sum_{\tau
_2=-L_{d_n}}^{L_{d_n}}
\sum_{\tau_1=1}^{L_{n}}\sum_{j_1=c_n+d_n-\tau_1+1}^{c_n+d_n} \varphi
_{j_1 j_2 \tau_1 \tau_2}, \\
d_3(n) & =& \frac{1}{d_n L_{d_n}}
\sum_{\tau_1=-L_{n}}^{-1}\sum_{j_1=c_n+1}^{c_n-\tau_1}
\sum_{\tau_2=-L_{n}}^{-1}\sum_{j_2=c_n+1}^{c_n-\tau_2} \varphi_{j_1
j_2 \tau_1 \tau_2}\\
&&{}
+ \frac{1}{d_n L_{d_n}} \sum_{\tau_1=-L_{n}}^{-1}\sum
_{j_1=c_n+1}^{c_n-\tau_1}
\sum_{\tau_2=1}^{L_{n}}\sum_{j_2=c_n+d_n-\tau_2+1}^{c_n+d_n} \varphi
_{j_1 j_2 \tau_1 \tau_2} \\
&&{}+  \frac{1}{d_n L_{d_n}} \sum_{\tau_1=1}^{L_{n}}\sum
_{j_1=c_n+d_n-\tau_1+1}^{c_n+d_n}
\sum_{\tau_2=-L_{n}}^{-1}\sum_{j_2=c_n+1}^{c_n-\tau_2} \varphi_{j_1
j_2 \tau_1 \tau_2}\\
&&{}
+ \frac{1}{d_n L_{d_n}} \sum_{\tau_1=1}^{L_{n}}\sum_{j_1=c_n+d_n-\tau
_1+1}^{c_n+d_n}
\sum_{\tau_2=1}^{L_{n}}\sum_{j_2=c_n+d_n-\tau_2+1}^{c_n+d_n} \varphi
_{j_1 j_2 \tau_1 \tau_2}
%
\end{eqnarray*}
and $\varphi_{j_1 j_2
\tau_1 \tau_2}= E(X_{j_1}X_{j_2})E(X_{j_1+\tau_1}X_{j_2+\tau_2})
H_{L_{d_n}}(\tau_1)H_{L_{d_n}}(\tau_2)
\mathrm{e}^{-\mathrm{i}(\nu_1 (j_1+\tau_1) - \omega_1 j_1 - \nu_2
(j_2+\tau_2) + \omega_2 j_2)}$.
We show that the terms $d_{1}(n), d_{2}(n), d_{3}(n)$ tend to zero,
which means that $s_{2}(n)$ has the same limit
as $d_{0}(n)$. We start with the term $d_{1}(n)$. By Lemma \ref
{base}(i), we have
$|\varphi_{j_1 j_2 \tau_1 \tau_2 }| \leq64 \Delta^4
\alpha^{\delta/{(2+\delta)}}(|j_1-j_2|)\alpha^{\delta
/{(2+\delta)}}(|j_1-j_2-(\tau_1-\tau_2)|)$,
which means that
\begin{eqnarray*}
|d_1(n)| & \leq& \frac{64 \Delta^4 }{d_n L_{d_n}} \sum_{j_1=c_n+1}^{c_n+d_n}
\sum_{\tau_1=-L_{d_n}}^{L_{d_n}}
\sum_{\tau_2=-L_{n}}^{-1}
\sum_{j_2=c_n+1}^{c_n-\tau_2}\alpha^{\delta/{(2+\delta)}}(|j_1-j_2|)\\[-1pt]
&&\hphantom{\frac{64 \Delta^4 }{d_n L_{d_n}} \sum_{j_1=c_n+1}^{c_n+d_n}
\sum_{\tau_1=-L_{d_n}}^{L_{d_n}}
\sum_{\tau_2=-L_{n}}^{-1}
\sum_{j_2=c_n+1}^{c_n-\tau_2}}
{}\times\alpha^{\delta
/{(2+\delta)}}\bigl(|j_1-j_2-(\tau_1-\tau_2)|\bigr)\\[-1pt]
&&{} + \frac{64 \Delta^4 }{d_n L_{d_n}} \sum_{j_1=c_n+1}^{c_n+d_n}
\sum_{\tau_1=-L_{d_n}}^{L_{d_n}}
\sum_{\tau_2=1}^{L_{n}}
\sum_{j_2=c_n+d_n-\tau_2+1}^{c_n+d_n}
\alpha^{\delta/{(2+\delta)}}(|j_1-j_2|)\\[-1pt]
&&\hphantom{{} + \frac{64 \Delta^4 }{d_n L_{d_n}} \sum_{j_1=c_n+1}^{c_n+d_n}
\sum_{\tau_1=-L_{d_n}}^{L_{d_n}}
\sum_{\tau_2=1}^{L_{n}}
\sum_{j_2=c_n+d_n-\tau_2+1}^{c_n+d_n}}
{}\times\alpha^{\delta
/{(2+\delta)}}\bigl(|j_1-j_2-(\tau_1-\tau_2)|\bigr)\\[-1pt]
& \leq& \frac{64 \Delta^4 }{d_n L_{d_n}} \Biggl( \sum_{k=1-L_{d_n}}^{d_n-1}
L_{d_n} \alpha^{\delta/{(2+\delta)}}(|k|)
\sum_{\tau_1=-L_{d_n}}^{L_{d_n}}
\sum_{\tau_2=-L_{d_n}}^{L_{d_n}} \alpha^{\delta/{(2+\delta)
}}\bigl(|k-(\tau_1-\tau_2)|\bigr) \\[-1pt]
&&\hphantom{\frac{64 \Delta^4 }{d_n L_{d_n}} \Biggl(}{} + \sum_{k=-d_n+1}^{L_n-1} L_{d_n} \alpha^{\delta/{(2+\delta)}}(|k|)
\sum_{\tau_1=-L_{d_n}}^{L_{d_n}}
\sum_{\tau_2=-L_{d_n}}^{L_{d_n}} \alpha^{\delta/{(2+\delta)
}}\bigl(|k-(\tau_1-\tau_2)|\bigr) \Biggr) \\[-1pt]
& \leq& \frac{64 \Delta^4 }{d_n L_{d_n}} \Biggl( \sum_{k=1-L_{d_n}}^{d_n-1}
L_{d_n} \alpha^{\delta/{(2+\delta)}}(|k|) (2
L_{d_n} +1)2 K\\[-1pt]
&&\hphantom{\frac{64 \Delta^4 }{d_n L_{d_n}} \Biggl(}
{} + \sum_{k=-d_n+1}^{L_n-1} L_{d_n} \alpha^{\delta
/{(2+\delta)}}(|k|) (2 L_{d_n}+1)2 K \Biggr)\\[-1pt]
& \leq&\frac{64 \Delta^4 }{d_n L_{d_n}} (2 L_{d_n}+1)^2 8 K^2 \to0
\qquad \mbox{as } n \to\infty.
\end{eqnarray*}
%
Using the same steps, we have that the terms $d_{2}(n)$ and $d_{3}(n)$
tend to zero. For $d_{0}(n)$, we get
\begin{eqnarray*}
d_{0}(n) & = &\frac{1}{d_n L_{d_n}} \sum_{j_1=c_n+1}^{c_n+d_n} \sum
_{j_2=c_n+1}^{c_n+d_n} E(X_{j_1}X_{j_2})\mathrm{e}^{-\mathrm{i}( \omega_2 j_2 - \omega_1 j_1 )}\\[-1pt]
&&\hphantom{\frac{1}{d_n L_{d_n}} \sum_{j_1=c_n+1}^{c_n+d_n} \sum
_{j_2=c_n+1}^{c_n+d_n}}
{}\times
\sum_{\tau_1=-L_{d_n}}^{L_{d_n}} \sum_{\tau_2=-L_{d_n}}^{L_{d_n}}
E(X_{j_1+\tau_1} X_{ j_1+\tau_2 })\\[-1pt]
&&\hphantom{\hphantom{\frac{1}{d_n L_{d_n}} \sum_{j_1=c_n+1}^{c_n+d_n} \sum
_{j_2=c_n+1}^{c_n+d_n}}
{}\times\sum_{\tau_1=-L_{d_n}}^{L_{d_n}} \sum_{\tau_2=-L_{n}}^{L_{d_n}}}
{}\times H_{L_{d_n}}(\tau_1)^2 \mathrm{e}^{-\mathrm{i}(\nu_1 (j_1+\tau_1) - \nu
_2 (j_1+\tau_2))}\\
&&{} + \frac{1}{d_n L_{d_n}} \sum_{j_1=c_n+1}^{c_n+d_n} \sum
_{j_2=c_n+1}^{c_n+d_n} E(X_{j_1}X_{j_2})
\mathrm{e}^{-\mathrm{i}( \omega_2 j_2 - \omega_1 j_1 )}\\
&&\hphantom{{} + \frac{1}{d_n L_{d_n}} \sum_{j_1=c_n+1}^{c_n+d_n} \sum
_{j_2=c_n+1}^{c_n+d_n}}
{}\times
\sum_{\tau_1=-L_{d_n}}^{L_{d_n}} \sum_{\tau_2=-L_{d_n}}^{L_{d_n}}
\bigl( E(X_{j_1+\tau_1} X_{ j_2+\tau_2 }) \mathrm{e}^{ \mathrm{i} \nu_2 (j_2+\tau_2))} \\
&&\hphantom{\hphantom{{} + \frac{1}{d_n L_{d_n}} \sum_{j_1=c_n+1}^{c_n+d_n} \sum
_{j_2=c_n+1}^{c_n+d_n}}
{}\times
\sum_{\tau_1=-L_{d_n}}^{L_{d_n}} \sum_{\tau_2=-L_{n}}^{L_{d_n}}
\bigl( }
{} - E(X_{j_1+\tau_1} X_{ j_1+\tau_2 }) \mathrm{e}^{ \mathrm{i} \nu_2 (j_1+\tau_2))}
\bigr)\\
&&\hphantom{\hphantom{{} + \frac{1}{d_n L_{d_n}} \sum_{j_1=c_n+1}^{c_n+d_n} \sum
_{j_2=c_n+1}^{c_n+d_n}}
{}\times
\sum_{\tau_1=-L_{d_n}}^{L_{d_n}} \sum_{\tau_2=-L_{n}}^{L_{d_n}}}
{}\times
H_{L_{d_n}}(\tau_1)^2
\mathrm{e}^{-\mathrm{i} \nu_1 (j_1+\tau_1) }\\
&&{} + \frac{1}{d_n L_{d_n}} \sum_{j_1=c_n+1}^{c_n+d_n} \sum
_{j_2=c_n+1}^{c_n+d_n} E(X_{j_1}X_{j_2})
\mathrm{e}^{-\mathrm{i}( \omega_2 j_2 - \omega_1 j_1 )}\\
&&\hphantom{{} + \frac{1}{d_n L_{d_n}} \sum_{j_1=c_n+1}^{c_n+d_n} \sum
_{j_2=c_n+1}^{c_n+d_n}}
{}\times\sum_{\tau_1=-L_{d_n}}^{L_{d_n}} \sum_{\tau_2=-L_{d_n}}^{L_{d_n}}
E(X_{j_1+\tau_1} X_{ j_2+\tau_2 })\\
&&\hphantom{\hphantom{{} + \frac{1}{d_n L_{d_n}} \sum_{j_1=c_n+1}^{c_n+d_n} \sum
_{j_2=c_n+1}^{c_n+d_n}}
{}\times\sum_{\tau_1=-L_{d_n}}^{L_{d_n}} \sum_{\tau_2=-L_{n}}^{L_{d_n}}}
{}\times H_{L_{d_n}}(\tau_1) \bigl( H_{L_{d_n}}(\tau_2) - H_{L_{d_n}}(\tau_1)
\bigr)\\
&&\hphantom{\hphantom{{} + \frac{1}{d_n L_{d_n}} \sum_{j_1=c_n+1}^{c_n+d_n} \sum
_{j_2=c_n+1}^{c_n+d_n}}
{}\times\sum_{\tau_1=-L_{d_n}}^{L_{d_n}} \sum_{\tau_2=-L_{n}}^{L_{d_n}}}
{}\times\mathrm{e}^{-\mathrm{i}(\nu_1 (j_1+\tau_1) - \nu_2 (j_2+\tau_2))}.
%
\end{eqnarray*}
Denote the first, second and third terms of the last equality by
$c_0(n)$, $c_1(n)$ and $c_2(n)$, respectively.
For the term $c_1(n)$, we have $c_1(n)=z_1(n)+z_2(n)+z_3(n)$, where
\begin{eqnarray*}
z_1 (n) & = &\frac{1}{d_n L_{d_n}}\sum_{|j_1-j_2|\geq2 L_{d_n}}
E(X_{j_1}X_{j_2})
\mathrm{e}^{-\mathrm{i}( \omega_2 j_2 - \omega_1 j_1 )}\\
&&\hphantom{\frac{1}{d_n L_{d_n}}\sum_{|j_1-j_2|\geq2 L_{d_n}}}
{}\times\sum_{\tau_1=-L_{d_n}}^{L_{d_n}} \sum_{\tau_2=-L_{d_n}}^{L_{d_n}}
\bigl( E(X_{j_1+\tau_1} X_{ j_2+\tau_2 }) \mathrm{e}^{ \mathrm{i} \nu_2 (j_2+\tau_2))} \\
&&\hphantom{\hphantom{\frac{1}{d_n L_{d_n}}\sum_{|j_1-j_2|\geq2 L_{d_n}}}
{}\times\sum_{\tau_1=-L_{d_n}}^{L_{d_n}} \sum_{\tau_2=-L_{n}}^{L_{d_n}}
\bigl(}
{} - E(X_{j_1+\tau_1} X_{ j_1+\tau_2 }) \mathrm{e}^{ \mathrm{i} \nu_2 (j_1+\tau_2))}
\bigr)\\
&&\hphantom{\hphantom{\frac{1}{d_n L_{d_n}}\sum_{|j_1-j_2|\geq2 L_{d_n}}}
{}\times\sum_{\tau_1=-L_{d_n}}^{L_{d_n}} \sum_{\tau_2=-L_{n}}^{L_{d_n}}}
{}\times
H_{L_{d_n}}(\tau_1)^2
\mathrm{e}^{-\mathrm{i} \nu_1 (j_1+\tau_1) },\\
z_2 (n) & =& \frac{1}{d_n L_{d_n}}\sum_{0 < j_1 - j_2 < 2
L_{d_n}}E(X_{j_1}X_{j_2})
\sum_{\tau_1=-L_{d_n}}^{L_{d_n}}
\Biggl( \sum_{\tau_2=j_2-j_1-L_{d_n}}^{-L_{d_n}-1} -
\sum_{\tau_2=j_2-j_1+L_{d_n}+1}^{L_{d_n}} \Biggr)\\
&&\hphantom{\frac{1}{d_n L_{d_n}}\sum_{0 < j_1 - j_2 < 2
L_{d_n}}E(X_{j_1}X_{j_2})
\sum_{\tau_1=-L_{n}}^{L_{d_n}}}
{}\times E(X_{j_1+\tau_1} X_{ j_1+\tau_2 })
H_{L_{d_n}}(\tau_1)^2 \\
&&\hphantom{\frac{1}{d_n L_{d_n}}\sum_{0 < j_1 - j_2 < 2
L_{d_n}}E(X_{j_1}X_{j_2})
\sum_{\tau_1=-L_{n}}^{L_{d_n}}}
{}\times\mathrm{e}^{-\mathrm{i}(\nu_1 (j_1+\tau_1) -
\omega_1 j_1 - \nu_2 (j_1+\tau_2) + \omega_2 j_2)},\\
z_3 (n) & =& \frac{1}{d_n L_{d_n}}\sum_{0 < j_2 - j_1 < 2
L_{d_n}}E(X_{j_1}X_{j_2})
\sum_{\tau_1=-L_{d_n}}^{L_{d_n}}
\Biggl(\sum_{\tau_2=L_{d_n}+1}^{j_2-j_1+L_{d_n}} -
\sum_{\tau_2=-L_{d_n}}^{j_2-j_1-L_{d_n}-1} \Biggr)\\
&&\hphantom{\frac{1}{d_n L_{d_n}}\sum_{0 < j_2 - j_1 < 2
L_{d_n}}E(X_{j_1}X_{j_2})
\sum_{\tau_1=-L_{n}}^{L_{d_n}}}
{}\times E(X_{j_1+\tau_1} X_{ j_1+\tau_2 })
H_{L_{d_n}}(\tau_1)^2 \\
&&\hphantom{\frac{1}{d_n L_{d_n}}\sum_{0 < j_2 - j_1 < 2
L_{d_n}}E(X_{j_1}X_{j_2})
\sum_{\tau_1=-L_{n}}^{L_{d_n}}}
{}\times \mathrm{e}^{-\mathrm{i}(\nu_1 (j_1+\tau_1) -
\omega_1 j_1 - \nu_2 (j_1+\tau_2) + \omega_2 j_2)}.
\end{eqnarray*}
Using Lemma \ref{base}(i) and similar steps as for the term $d_1 (n)$,
it can be proven that $z_{1}(n)$ tends to
zero. For the term $z_{2}(n)$, we have 
\begin{eqnarray*}
|z_2 (n)| & \leq&\frac{64\Delta^4}{d_n L_{d_n}}\sum_{0 < j_1 - j_2 <
2 L_{d_n}}
\sum_{\tau_1=-L_{d_n}}^{L_{d_n}}
\Biggl( \sum_{\tau_2=j_2-j_1-L_{d_n}}^{-L_{d_n}-1} +
\sum_{\tau_2=j_2-j_1+L_{d_n}+1}^{L_{d_n}} \Biggr)\\
&&\hphantom{\frac{64\Delta^4}{d_n L_{d_n}}\sum_{0 < j_1 - j_2 <
2 L_{d_n}}
\sum_{\tau_1=-L_{n}}^{L_{d_n}}}{}\times\alpha^{\delta/{(2+\delta)}}(j_1 - j_2) \alpha^{\delta
{(2+\delta)}}(|\tau_1 - \tau_2|)\\
& \leq&\frac{128K\Delta^4}{d_n L_{d_n}}\sum_{0 < j_1 - j_2 < 2
L_{d_n}}(j_1-j_2)
\alpha^{\delta/{(2+\delta)}}(j_1 - j_2) \to0.
\end{eqnarray*}
In the same way, we may prove that $z_3 (n)\to0$. This means that
$c_{1}(n) \to0$. Using Lemma \ref{base}(i),
inequality $|H_{L_{d_n}}(\tau_2) - H_{L_{d_n}}(\tau_1)|\leq W |\tau
_2-\tau_1|/L_{d_n} \leq W (|\tau_2-\tau_1 +
(j_2-j_1)|+|j_2-j_1|)/L_{d_n}$ and similar steps as for the term $d_1
(n)$, it can be proven that $c_{2}(n)$
tends to zero. This means that the term $d_{0}(n)$ has the same limit
as $c_{0}(n)$. Now, using Lemmas \ref{bias}
and \ref{biasx} for the term $c_{0}(n)$, we get
\begin{eqnarray}
c_{0}(n) & =& 2\uppi\int_{-1}^{1}w^2(x)\,\mathrm{d}x\, P(\nu_1,\nu_2) \frac{1}{d_n}
\sum_{j_1=c_n+1}^{c_n+d_n} \sum_{j_2=c_n+1}^{c_n+d_n} B(j_1,j_2-j_1)
\mathrm{e}^{-\mathrm{i}(- \omega_1 j_1 + \omega_2 j_2)} \nonumber\\
&&{} + \underbrace{\frac{2\uppi}{d_n}
\sum_{j_1=c_n+1}^{c_n+d_n} \sum_{j_2=c_n+1}^{c_n+d_n} B(j_1,j_2-j_1)
\mathrm{e}^{-\mathrm{i}(- \omega_1 j_1 + \omega_2 j_2)} e_{n}(\nu_1,\nu
_2)}_{y_{n}} \\
%
& = &(2\uppi)^2 \int_{-1}^{1}w^2(x)\,\mathrm{d}x\, P(\nu_1,\nu_2)\overline
{P(\omega_1,\omega_2)} + \mathrm{o}(1) + y_{n}.\nonumber
\end{eqnarray}
For the term $y_{n}$, we have
\[
|y_{n}| \leq\tilde e_{n}(\nu_1,\nu_2) \frac{2\uppi}{d_n}
\sum_{j_1=c_n+1}^{c_n+d_n} \sum_{j_2=c_n+1}^{c_n+d_n} 8 \Delta^2
\alpha^{\delta/{(2+\delta)}}(|j_1-j_2|)
\leq\tilde e_{n}(\nu_1,\nu_2) 4 K \uppi\to0.
\]
Hence, $s_{2}(n)\to(2\uppi)^2 \rho P(\nu_1,\nu_2) \overline{P(\omega
_1,\omega_2)}$. Following the same steps,
we get that $s_{3}(n) \to (2\uppi)^2 \rho P(\nu_1,2\uppi-
\omega_2) \overline{P(\nu_2, 2\uppi- \omega_1)}$. This completes the proof.
\end{pf*}

\begin{pf*}{Proof of Theorem \ref{normality}}
Let us consider the following decomposition for the estimator $\hat
G_{n}^{c,d} (\nu,\omega)$:
\begin{eqnarray*}
\sqrt{d_n/L_{d_n}}\bigl(\hat G_{n}^{c,d}(\nu,\omega)- P(\nu,\omega)\bigr)&=&
\sqrt{d_n/L_{d_n}} \bigl(\hat
G_{n}^{c,d}(\nu,\omega)- E(\hat G_{n}^{c,d}(\nu,\omega)) \bigr)\\
&&{} +\sqrt{d_n/L_{d_n}}\bigl(E(\hat G_{n}^{c,d}(\nu,\omega))- P(\nu,\omega)\bigr)\\
&=& S^{c,d}_{n}(\nu,\omega)+ \epsilon^{c,d}_{n}(\nu,\omega),
\end{eqnarray*}
where $S^{c,d}_{n}(\nu,\omega) = \sqrt{d_n/L_{d_n}} (\hat
G_{n}^{c,d}(\nu,\omega)-
E(\hat G_{n}^{c,d}(\nu,\omega)) )$,
$\epsilon^{c,d}_{n}(\nu,\omega)= \sqrt{d_n/L_{d_n}}(E(\hat
G_{n}^{c,d}(\nu,\break\omega))- P(\nu,\omega)).$
We will now split the proof into two steps. In the first step, we will
show that the deterministic
term $\epsilon^{c,d}_{n}(\nu,\omega)$ tends to zero for $n\rightarrow
\infty$.
In the second step, we will show that
\begin{equation}\label{S_as_nor}
\left[
\matrix{
\operatorname{Re}(S^{c,d}_{n}(\nu,\omega))\cr
\operatorname{Im}(S^{c,d}_{n}(\nu,\omega))
}
\right]
\stackrel{d}{\longrightarrow}\mathcal{N}_{2}(0,\Sigma(\nu,\omega)).
\end{equation}

\textit{Step} 1. Without loss of generality, we may assume that $\omega
\leq \nu$.
Changing variables and using, sequentially, (\ref{fourier2}), (\ref
{spec_rec}), (\ref{as-per}), (A3) and assumption (ii) of our theorem,
we then obtain
\begin{eqnarray*}
|\epsilon_{n}^{c,d}(\nu,\omega)| &=& \sqrt{d_n/L_{d_n}} |E(\hat
G_{n}^{c,d}(\nu,\omega))-P(\nu,\omega)|\\
& = &\sqrt{\frac{d_n}{L_{d_n}}} \Biggl| \frac{1}{2\uppi d_n}\sum_{j=c_n+1}^{c_n+d_n}
\sum_{\tau=-L_{d_n}}^{L_{d_n}} \mathbf{1}\{c_n + 1 \leq j + \tau
\leq c_n + d_n \} B(j,\tau) H_{L_{d_n}}(\tau)\\
&&\hphantom{\sqrt{\frac{d_n}{L_{d_n}}} \Biggl| \frac{1}{2\uppi d_n}\sum_{j=c_n+1}^{c_n+d_n}
\sum_{\tau=-L_{d_n}}^{L_{d_n}}}
{}\times \mathrm{e}^{-\mathrm{i}\nu \tau}\mathrm{e}^{- \mathrm{i} (\nu - \omega) j} -P(\nu,\omega) \Biggr|\\
& \leq& \sqrt{\frac{d_n}{L_{d_n}}} \Biggl| \frac{1}{2\uppi d_n}\sum
_{j=c_n+1}^{c_n+d_n}
\sum_{\tau=-L_{d_n}}^{L_{d_n}} B(j,\tau) H_{L_{d_n}}(\tau) \mathrm
{e}^{-\mathrm{i} \nu \tau}\mathrm{e}^{- \mathrm{i} (\nu - \omega) j} -P(\nu,\omega) \Biggr|\\
&&{} + \sqrt{\frac{d_n}{L_{d_n}}} \frac{1}{2\uppi d_n}\sum_{j=c_n+1}^{c_n+d_n}
\sum_{\tau=-L_{d_n}}^{L_{d_n}} (\mathbf{1}\{ j+\tau> c_n + d_n \}\\
&&\hphantom{{} + \sqrt{\frac{d_n}{L_{d_n}}} \frac{1}{2\uppi d_n}\sum_{j=c_n+1}^{c_n+d_n}
\sum_{\tau=-L_{d_n}}^{L_{d_n}} (}
{}+\mathbf{1}\{ j+\tau< c_n+1 \}) 8 \Delta^2 \alpha^{\delta/{(2+\delta)}}(|\tau|)\\
& \leq& \sqrt{\frac{d_n}{L_{d_n}}} \biggl| \frac{1}{2\uppi}
\sum_{|\tau| \leq L_{d_n}} a(\nu-\omega,\tau) H_{L_{d_n}}(\tau)
\mathrm{e}^{-\mathrm{i} \nu\tau}
- P(\nu,\omega) \biggr|\\
&&{} + \sqrt{\frac{d_n}{L_{d_n}}} \Biggl| \frac{1}{2\uppi d_n} \sum_{\tau
=-L_{d_n}}^{L_{d_n}}
H_{L_{d_n}}(\tau) \sum_{j=c_n+1}^{c_n+d_n}
\sum_{\lambda\in\Lambda_\tau\setminus\{ \nu-\omega\}}
a(\lambda,\tau) \mathrm{e}^{\mathrm{i}(\lambda-(\nu-\omega))j} \mathrm
{e}^{-\mathrm{i} \nu\tau} \Biggr|\\
&&{} + \frac{4 \Delta^2}{\uppi\sqrt{d_n
L_{d_n}}} \sum_{\tau=-L_{d_n}}^{L_{d_n}} |\tau|\alpha^{\delta
/{(2 +\delta)}}(|\tau|)\\
& \leq &\sqrt{\frac{d_n}{L_{d_n}}} \biggl| \frac{1}{2\uppi} \biggl(
\sum_{|\tau| \leq\theta L_{d_n}} + \sum_{L_{d_n} >|\tau| > \theta
L_{d_n}} \biggr) a(\nu-\omega,\tau) H_{L_{d_n}}(\tau) \mathrm{e}^{-\mathrm{i}
\nu\tau}
- P(\nu,\omega) \biggr|\\
&&{} + \sqrt{\frac{d_n}{L_{d_n}}} \frac{1}{2\uppi d_n} \sum_{\tau
=-L_{d_n}}^{L_{d_n}} \sum_{\lambda\in\Lambda_\tau\setminus\{ \nu
-\omega\}}
\Biggl| \sum_{j=c_n+1}^{c_n+d_n}
a(\lambda,\tau) \mathrm{e}^{\mathrm{i}(\lambda-(\nu-\omega))j} \Biggr| +
\mathrm{o}(1) \\
& \leq& \sqrt{\frac{d_n}{L_{d_n}}} \biggl( \frac{1}{\uppi}
\sum_{|\tau|> \theta L_{d_n}} 8 \Delta^2 \alpha^{\delta
/{(2+\delta)}}(|\tau|)\\
&&\hphantom{\sqrt{\frac{d_n}{L_{d_n}}} \biggl(}
{} +
\frac{1}{2\uppi d_n} \sum_{\tau=-L_{d_n}}^{L_{d_n}}
\sum_{\lambda\in\Lambda_\tau\setminus\{ \nu-\omega\}}
\bigl|a(\lambda,\tau) \cosec\bigl(\bigl(\lambda- (\nu - \omega)\bigr)/2\bigr)\bigr| \biggr) + \mathrm{o}(1)\\
& \leq&\frac{8 \Delta^2}{ \uppi} \sqrt{\frac{d_n}{L_{d_n}}}
\sum_{|\tau|> \theta L_{d_n}} \alpha^{\delta/{(r+2+\delta)}}(|\tau
|)\\
&&{} +
\frac{1}{2\uppi\sqrt{d_n L_{d_n}}} \sum_{\tau=-L_{d_n}}^{L_{d_n}}
\sum_{\lambda\in\Lambda_\tau\setminus\{ \nu-\omega\}}
\bigl|a(\lambda,\tau) \cosec\bigl(\bigl(\lambda- (\nu - \omega)\bigr)/2\bigr)\bigr| + \mathrm{o}(1)\\
& \leq&\frac{8 \Delta^2}{\uppi\theta^r} \sqrt{\frac{d_n}{L_{d_n}^{1+2r}}}
\sum_{|\tau|>\theta L_{d_n}} \theta^r L_{d_n}^{r} \alpha^{\delta
/{(r+2+\delta)}}(|\tau|) \\
&&{}+
\frac{1}{2\uppi} \sqrt{\frac{L_{d_n}}{d_n}} \Biggl(\frac{1}{L_{d_n}}
\sum_{\tau=-L_{d_n}}^{L_{d_n}} z_{\tau}(\nu-\omega) \Biggr) + \mathrm{o}(1)
\\
& \leq&\frac{8 \Delta^2}{ \uppi\theta^r} \sqrt{\frac
{d_n}{L_{d_n}^{1+2r}}} K + \frac{1}{2\uppi} \sqrt{\frac{L_{d_n}}{d_n}}
\mathrm{o}(1) + \mathrm{o}(1)= \mathrm{O}(1) d_n^{(1-\kappa(1+2r))/2} + \mathrm{o}(1) \\
&= &\mathrm{O}(1) d_n^{\kappa
((1-\kappa)/{(2\kappa)}-r )} + \mathrm{o}(1)
\to0,
\end{eqnarray*}
where $\mathbf{1}{\{B\}}$ is an indicator function of the event $B$.
This completes the proof of step 1.

\textit{Step} 2. In this step, we use Theorem \ref{covariance} and Lemma
\ref{as-norm} to show (\ref{S_as_nor}). Note that
\[
S^{c,d}_{ n}(\nu,\omega)= \frac{1}{2\uppi\sqrt{d_n}} \sum
_{j=1}^{d_n} \bigl(R_{n,j}(\nu,\omega)
+\im I_{n,j}(\nu,\omega) \bigr),
\]
where $R_{n,j}(\nu,\omega)=R_{n,j}$ and $I_{j}(\nu,\omega)=I_{n,j}$
are defined via
\begin{eqnarray*}
R_{n,j}(\nu,\omega)&=&\frac{1}{\sqrt{ L_{d_n}}}\sum_{\tau
=-L_{d_n}}^{L_{d_n}} H_{L_{d_n}}(\tau)
\mathbf{1}\{1\leq j+\tau\leq d_n \} Z_{j+c_n}(\tau)\cos\bigl(\nu \tau+
(\nu- \omega) (j+c_n)\bigr),
\\
I_{n,j}(\nu,\omega)&=&\frac{1}{\sqrt{ L_{d_n}}}\sum_{\tau
=-L_{d_n}}^{L_{d_n}} H_{L_{d_n}}(\tau)
\mathbf{1}\{1\leq j+\tau\leq d_n \} Z_{j+c_n}(\tau)\sin\bigl(- \nu \tau
- (\nu- \omega ) (j+c_n)\bigr)
\end{eqnarray*}
and
$Z_{j}(\tau)=X_{j}X_{j+\tau}-B(j,\tau)$. We show that for triangular
array $\{(R_{n,t},I_{n,t})\dvtx 1\leq t\leq d_n\}$, the assumptions of
Theorem \ref{as-norm} hold.
Note that by the Minkowski inequality, Lemma \ref{base}(i) and,
finally, by the H\"older inequality, we have
\begin{eqnarray*}\label{est1}
&&\| R_{n,j}(\nu,\omega)\|_{2+\delta} \\
&&\quad \leq \Biggl\|\frac{1}{\sqrt{ L_{d_n}}}
\sum_{\tau=-L_{d_n}}^{L_{d_n}} H_{L_{d_n}}(\tau)
\mathbf{1}\{1\leq j + \tau\leq d_n \}
X_{j + c_n}X_{j + c_n + \tau} \cos\bigl(\nu \tau + (\nu - \omega) (j +
c_n)\bigr) \Biggr\|_{2+\delta}\\
& &\qquad {}+ \frac{8 \Delta^2}{\sqrt{ L_{d_n}}} \sum_{\tau=-L_{d_n}}^{L_{d_n}}
\alpha^{\delta/{(2+\delta)}}(|\tau|) \\
& &\quad  \leq\| X_{j + c_n} \|_{6+3\delta} \Biggl\|\frac{1}{\sqrt{
L_{d_n}}}\sum_{\tau=-L_{d_n}}^{L_{d_n}}
H_{L_{d_n}}(\tau)
\mathbf{1}\{1\leq j + \tau\leq d_n \}
X_{j + c_n + \tau} \\
&&\hphantom{\quad  \leq\| X_{j + c_n} \|_{6+3\delta} \Biggl\|\frac{1}{\sqrt{L_{d_n}}}\sum_{\tau=-L_{d_n}}^{L_{d_n}}}
{}\times\cos\bigl(\nu \tau + (\nu - \omega) (j + c_n)\bigr) \Biggr\|_{3
+ 3 \delta/2 }
+ 16 \Delta^2 K \\
& &\quad  \leq\frac{\Delta}{\sqrt{ L_{d_n}}}
\underbrace{ \Biggl\|\sum_{\tau=-L_{d_n}}^{L_{d_n}}
H_{L_{d_n}}(\tau)
\mathbf{1}\{1\leq j+\tau\leq d_n \}
X_{j+c_n+\tau}
\cos\bigl(\nu \tau+ (\nu- \omega) (j+c_n)\bigr) \Biggr\|_{3 + 3 \delta/2 }
}_{t_1(n)}\\
&&\qquad {} + 16 \Delta^2 K.
%
%
\end{eqnarray*}
In the next step, we use Lemma \ref{fazekas} to estimate the term
$t_1(n)$. Let $l=3 + \frac{3}{2}\delta$ and note that $l+\delta
<6+3\delta$, which means that assumption (i) of Lemma~\ref{fazekas} holds.
Assumption (ii) of Lemma~\ref{fazekas} follows from assumption (iii) of our theorem. Therefore,
using Lemma~\ref{fazekas} (similarly as for $I_{n,j}(\nu,\omega)$),
we get
\begin{eqnarray*}\label{est1}
\hspace*{-4pt}&&\| R_{n,j}(\nu,\omega)\|_{2+\delta}\\[2pt]
\hspace*{-4pt}&&\quad  \leq \Delta\bigl(K_{\alpha,l}
\max\bigl\{ Q(3 + 3\delta/2,\delta,-L_{d_n}-1,2 L_{d_n}+1),\\[2pt]
\hspace*{-4pt}&&\qquad\hphantom{\Delta\bigl(K_{\alpha,l}\max\{}
[Q(2,\delta,-L_{d_n}-1,2 L_{d_n}+1)]^{(3 + (3/2)\delta)/2} \bigr\}
\bigr)^{1/(3 + (3/2)\delta)}/\sqrt{L_{d_n}} + 16 \Delta^2 K \\[2pt]
\hspace*{-4pt}&&\quad  \leq\Delta K_{\alpha,l}^{1/(3 + (3/2)\delta)} \bigl( \max\bigl\{
(2 L_{d_n}+1) \Delta^{3 + (3/2)\delta},\bigl((2 L_{d_n}+1) \Delta
^{2}\bigr)^{(3 + (3/2)\delta)/2} \bigr\}\bigr)^{1/(3 + (3/2)\delta)}/\sqrt
{L_{d_n}}\\[2pt]
\hspace*{-4pt}&&\qquad {} + 16 \Delta^2 K \\[2pt]
\hspace*{-4pt}&&\quad  \leq\Delta^2 K_{\alpha,l}^{1/(3 + (3/2)\delta)} \sqrt{2
L_{d_n}+1} /\sqrt{L_{d_n}} +
16 \Delta^2 K \\[2pt]
\hspace*{-4pt}&&\quad\leq\Delta^2 \bigl(2 K_{\alpha,l}^{1/(3 + (3/2)\delta
)} + 16 K\bigr),
%
%
\end{eqnarray*}
which means that assumption (i) of Lemma \ref{as-norm} holds. Condition
(ii) of Lemma \ref{as-norm} follows
from Theorem \ref{covariance} by noting that
\begin{eqnarray*}
&&\cov(\operatorname{Re}(S^{c,d}_{n}(\nu,\omega)),\operatorname
{Re}(S^{c,d}_{n}(\nu,\omega)))\\
&&\quad = \frac{d_n}{4 L_{d_n}} \cov \bigl(\hat G_{n}^{c,d}(\nu,\omega) + \overline
{\hat G_{n}^{c,d}(\nu,\omega)},\hat G_{n}^{c,d}(\nu,\omega)
+ \overline{\hat G_{n}^{c,d}(\nu,\omega)} \bigr),\\
&&\cov (\operatorname{Re}(S^{c,d}_{n}(\nu,\omega)),\operatorname{Im}(S^{c,d}_{n}(\nu,\omega))) \\
&&\quad = \frac{d_n}{4 \im L_{d_n}} \cov\bigl(\hat G_{n}^{c,d}(\nu,\omega) +
\overline{\hat G_{n}^{c,d}(\nu,\omega)},\hat G_{n}^{c,d}(\nu,\omega)
- \overline{\hat G_{n}^{c,d}(\nu,\omega)} \bigr),\\
&&\cov (\operatorname{Im}(S^{c,d}_{n}(\nu,\omega)),\operatorname
{Im}(S^{c,d}_{n}(\nu,\omega))) \\
&&\quad = \frac{d_n}{4 L_{d_n}} \cov\bigl(\hat G_{n}^{c,d}(\nu,\omega) - \overline
{\hat G_{n}^{c,d}(\nu,\omega)},\hat G_{n}^{c,d}(\nu,\omega)
- \overline{\hat G_{n}^{c,d}(\nu,\omega)} \bigr).
\end{eqnarray*}
Finally, putting $h_n=2 L_{d_n}$ and
$\varsigma(k)=\alpha_{X}(k)$, and taking into consideration assumption
(iii) of our theorem, we get that condition (iii)
of Lemma \ref{as-norm} holds,
which means that (\ref{S_as_nor}) holds,
where $\Sigma(\nu,\omega)$ can be obtained by the last three equations and
Theorem \ref{covariance}. The calculation of the matrix $\Sigma(\nu
,\omega)$ is too technical to be presented here.
This completes the proof.
\end{pf*}

\begin{pf*}{Proof of Theorem \ref{subsamplingp}}
By Politis \textit{et al.} \cite{Pol}, Theorem 4.2.1, it is
sufficient to prove that there exists a continuous distribution $J$
such that:
\begin{longlist}[(ii)]
\item[(i)] $\sqrt{n/L_{n}}(|\hat G_{n}(\nu,\omega)|-|P(\nu,\omega
)|)\stackrel{d}{\longrightarrow} J $;
\item[(ii)] for any sequence of positive
integers $\{t_{b}\}=\{t_{b(n)}\}$ such that $b=b(n) \to\infty$, $b/n
\to0$ as $n \to\infty$, we have
$J_{b,t_b}^{P}(x)\longrightarrow J(x),$
where
$J_{b,t_b}^{P}(x)=P(\sqrt{b/L_{b}}(|\hat G_{n}^{t_b-1,b}(\nu,\omega
)|-|P(\nu,\omega)|)\leq x)$
and $J(x)$ is a distribution function at the point $x$ for the
distribution $J$.
\end{longlist}
This follows immediately from Corollary \ref{J1-dis} by setting
$c_n=t_{b(n)}$ and $d_n=b(n)$. This completes the proof.
\end{pf*}

\begin{pf*}{Proof of Theorem \ref{subsamplingc}}
This proof is analogous to the proof of Theorem 4.1.
The only difference is that we use Corollary \ref{J2-dis-a} instead of
Corollary \ref{J1-dis}.

\end{pf*}

\end{appendix}

\section*{Acknowledgements}
My thanks go to
Dominique Dehay, Jacek Le\'skow and Rafa\l{} Synowiecki for
stimulating discussions. This research was supported in part by NATO
Grant ICS.NUKR.CLG 983335.

\printhistory

\end{document}